\documentclass[10pt]{amsart}
\usepackage[left=1in,top=1in,right=1in,bottom=1in]{geometry}
\usepackage{amsfonts}
\usepackage{enumerate,subfigure,color}
\usepackage{amssymb,latexsym, amsmath,graphicx}
\usepackage{epstopdf}
\usepackage{framed}
\usepackage{multirow}
\usepackage{booktabs}
\usepackage[round,sort&compress]{natbib}

\usepackage{hyperref}

\usepackage{enumitem}
\usepackage[utf8]{inputenc}
\usepackage{tikz}
\usetikzlibrary{shapes.geometric, arrows}

\tikzstyle{title} =  [rectangle, minimum width=3cm, minimum height=1cm,text centered, text width=5cm, draw=black, fill=purple!20]

\tikzstyle{startstop} = [rectangle, rounded corners, minimum width=3cm, minimum height=1cm,text centered, text width=5cm, draw=black, fill=orange!10]
\tikzstyle{io} =  [rectangle, rounded corners, minimum width=3cm, minimum height=1cm,text centered, text width=5cm, draw=black, fill=blue!20]
\tikzstyle{process} = [rectangle, rounded corners, minimum width=3cm, minimum height=1cm, text width=5cm, text centered, draw=black, fill=yellow!30]
\tikzstyle{algend} = [rectangle, rounded corners, minimum width=3cm, minimum height=1cm,text centered, text width=4.5cm, draw=black, fill=green!10]

\tikzstyle{decision} = [diamond, minimum width=3cm, minimum height=1cm, text centered, text width=3cm,draw=black, fill=green!30]

\tikzstyle{arrow} = [thick,->,>=stealth]

\tikzstyle{startstop2} = [rectangle, rounded corners, minimum width=3cm, minimum height=1cm,text centered, text width=8cm, draw=black, fill=orange!10]
\tikzstyle{io2} =  [rectangle, rounded corners, minimum width=3cm, minimum height=1cm,text centered, text width=7.5cm, draw=black, fill=blue!20]
\tikzstyle{process2} = [rectangle, rounded corners, minimum width=3cm, minimum height=1cm, text width=6cm, text centered, draw=black, fill=yellow!30]
\tikzstyle{algend2} = [rectangle, rounded corners, minimum width=3cm, minimum height=1cm,text centered, text width=10cm, draw=black, fill=green!10]

\tikzstyle{algend3} = [rectangle, rounded corners, minimum width=3cm, minimum height=1cm,text centered, text width=12.5cm, draw=black, fill=green!10]

\newcommand{\dbar}{\bar{\partial}}

\DeclareMathOperator{\T}{{\mathbf{t}}}

\DeclareMathOperator{\R}{{\mathbb R}}

\DeclareMathOperator{\C}{{\mathbb C}}

\DeclareMathOperator{\de}{\partial}

\DeclareMathOperator{\dez}{\de_z}
\DeclareMathOperator{\dbarz}{\dbar_z}
\DeclareMathOperator{\dbark}{\dbar_k}

\DeclareMathOperator{\by}{\times}
\DeclareMathOperator{\bndry}{\partial\Omega}

\DeclareMathOperator{\texp}{\mathbf{t}^{\mbox{\tiny{\textbf{exp}}}}}

\DeclareMathOperator{\sigDBn}{\sigma_n^{\mbox{\tiny \bf DB}}}

\newcommand{\D}{{\mathbb D}}

\newcommand{\tredbf}[1]{\textcolor{red}{\textbf{#1}}}

\usepackage{soul,color}
\soulregister\cite7
\soulregister\ref7
\soulregister\pageref7

\begin{document}
\bibliographystyle{plainnat} 

\title[Beltrami-Net for 2D Absolute EIT]{Beltrami-Net: Domain Independent Deep D-bar Learning for Absolute Imaging with Electrical Impedance Tomography (a-EIT)}

\author{S.~J. Hamilton,  A. H\"anninen, A. Hauptmann, and V. Kolehmainen}


\thanks{S.~J. Hamilton is with the Department of Mathematics, Statistics, and Computer Science; Marquette University, Milwaukee, WI 53233 USA,  email: \texttt{sarah.hamilton@marquette.edu}}
\thanks{A. H\"anninen is with the Department of Applied Physics; University of Eastern Finland, Kuopio, Finland, email: \texttt{asko.hanninen@uef.fi}}
\thanks{A. Hauptmann is with the Department of Computer Science; University College London, London, United Kingdom, email: \texttt{a.hauptmann@ucl.ac.uk}}
\thanks{V. Kohlehmainen is with the Department of Applied Physics; University of Eastern Finland, Kuopio, Finland, email: \texttt{ville.kolehmainen@uef.fi}}
\begin{abstract}
{\it Objective:}  To develop, and demonstrate the feasibility of, a novel image reconstruction method for absolute Electrical Impedance Tomography (a-EIT) that pairs deep learning techniques with real-time robust D-bar methods.  {\it Approach:}  A D-bar method is paired with a trained Convolutional Neural Network (CNN) as a post-processing step.  Training data is simulated for the network using no knowledge of the boundary shape by using an associated nonphysical Beltrami equation rather than simulating the traditional current and voltage data specific to a given domain. This allows the training  data to be boundary shape independent.  The method is tested on experimental data from two EIT systems (ACT4 and KIT4).  {\it Main Results:} Post processing the D-bar images with a CNN produces significant improvements in image quality measured by Structural SIMilarity indices (SSIMs) as well as relative $\ell_2$ and $\ell_1$ image errors.  {\it Significance:} This work demonstrates that more general networks can be trained without being specific about boundary shape, a key challenge in EIT image reconstruction.  The work is promising for future studies involving databases of anatomical atlases. 
\end{abstract}
\maketitle 


\section{Introduction}\label{sec:intro}

Electrical Impedance Tomography (EIT) probes a body with low-amplitude electrical currents applied on surface electrodes.  The surface measurements can then be used as inputs to solve a mathematical inverse problem to recover the internal electrical properties (conductivity and permittivity) of the object.  As EIT is a low-cost, non-invasive imaging modality with no ionizing radiation, it has several medical and industrial applications, see [\cite{Cheney1999}] and [\cite{Mueller2012}].  The image recovery task in EIT, recovering the internal conductivity from the surface electrode measurements, is a severely ill-posed nonlinear inverse problem thus requiring carefully designed reconstruction algorithms capable of handling incorrectly known boundary shape, electrode locations, and noise in the measured EIT data.  The ill-posedness of the inverse problem often results in images with low spatial resolution or severe image corruption due to modeling errors in a minimization task.  The D-bar method [\cite{Nachman1996, Knudsen2009}] has been shown to be robust to modeling errors and noise [\cite{Murphy2009}; \cite{Hamilton2018_Robust}].   

By viewing these low-resolution, real-time [\cite{Dodd2014}], D-bar images as convolutions of the true images one can develop and train a Convolutional Neural Network (CNN) to learn the blurring inherent in the D-bar reconstruction process on data of that type.  This idea was introduced in [\cite{Hamilton2018_DeepDbar}] and tested on experimental EIT data for absolute imaging in 2D.  There, the training data for the network was simulated from the forward EIT model
\begin{equation}\label{eq:eit_model}
\begin{array}{rcl}
\nabla\cdot \sigma(z)\nabla u(z) & = & 0, \quad z\in\Omega\subset\R^2\\
\sigma\frac{\partial u}{\partial \nu}&=& g, \quad z\in\bndry
\end{array}
\end{equation}
using the electrode continuum model [\cite{Hyvoenen2009}; \cite{Hauptmann2017}] based on continuum current/voltage data computed from a known circular domain boundary.  The trained network was then directly applied to D-bar reconstructions from the experimental data with no transfer training required.  By contrast, here we simulate our training data from the associated, non-physical, Beltrami problem [\cite{Astala2006a,Astala2006}] and `Shortcut D-bar Method' [\cite{Astala2010}] to remove any knowledge of the boundary (shape and electrodes) from the training process.  We test the network on EIT data from two different EIT machines (ACT4 [\cite{Liu2005}] and KIT4 [\cite{Kourunen2008}]) with different boundary shapes.  In practice, a network could be constructed using a database of CT scans where all that is needed is approximate internal structure boundaries (heart, lungs, spine, etc) and reasonable conductivity value windows for each type of inclusion. The CTs could be scaled such that the maximum radial component of the thorax boundary is one.  Alternatively, one could bypass any direct incorporation of organs by instead training using inclusions of ellipses, circles, etc.  The patient-specific voltage and current EIT data would then be scaled to correspond to a maximum radius of~1 by scaling the associated DN (or ND) matrix by the largest radial component of the patient's approximated boundary shape (see [\cite{Isaacson2004}]). In this study we investigate the particular question of how informative the training data needs to be in order to perform the desired image enhancement task after an initial reconstruction. That means, we consider two different scenarios in this study:
\begin{itemize}
\item[i.)] Thoracic measurements for a human patient, here a database can be built from anatomical atlases. In this setting the imaging task is highly constrained by anatomical features and hence training data can be tuned to be specific for this particular task. This constitutes a case of high \emph{a priori} knowledge.  We consider tank data with thoracic specific agar targets.
\item[ii.)] Assessment of the most general training data without any anatomical knowledge, with which we are able to achieve sufficient reconstruction quality for a vast application area. This can be considered a task of minimizing necessary \emph{a priori} knowledge.
\end{itemize}

The application of deep learning methods, in particular Convolutional Neural Networks (CNNs), has attracted major attention in recent years and shows great promise for improving images in tomographic reconstruction tasks. The most prominent approach is given by post-processing of an initial reconstruction based on an analytic inversion formula, such as filtered back-projection in X-ray CT [\cite{Kang2017}] and [\cite{Jin2017}].  Another promising clinical applications of this approach is dynamic cardiovascular magnetic resonance imaging [\cite{Schlemper2018}; \cite{Hauptmann2018}]. Recent studies, in addition to [\cite{Hamilton2018_DeepDbar}], have explored the possibility of using deep learning for EIT with artificial neural networks [\cite{Martin2017}] and variational autoencoders for lung imaging [\cite{Seo2018}].
Furthermore, several studies propose combining iterative variational techniques with deep learning to obtain superior reconstruction quality and more flexible generalization by including the forward operator in the network architectures [\cite{Adler2017}], [\cite{Hammernik2018}] and [\cite{Hauptmann2018a}].
\\

Section~\ref{sec:methods} presents the methods used in this work including the proposed new algorithm and how reconstruction quality will be assessed.  
Results of the proposed method on experimental EIT tank data from ACT4 and KIT4 are presented in Section~\ref{sec:results} and conclusions drawn in Section~\ref{sec:conc}.

\section{Methods}\label{sec:methods}
Here we consider the 2D real-valued conductivity EIT problem
\begin{equation}\label{eq:cond}
\nabla\cdot\sigma(z)\nabla u(z)=0,\qquad z\in\Omega\subset\R^2,
\end{equation}
where $\sigma=\sigma(z)$ is the spatially dependent conductivity and $u=u(z)$ the electric potential.  The current and voltage measurements take the form of approximate knowledge of the Neumann-to-Dirichlet (ND) map $\mathcal{R}_\sigma:\sigma\frac{\partial u}{\partial \nu} \mapsto g$ for $z\in\bndry$ which maps a boundary current to the corresponding boundary voltage, and $\nu=\nu(z)$ denotes the outward unit normal vector to $\bndry$.  Here, for simplicity, we assume the conductivity is constant $\sigma=\sigma_0$ in a neighborhood of the boundary.  If $\sigma $ is not constant near $\bndry$, a padding of the domain can be used as in [\cite{Nachman1996,Siltanen2014}] reducing the problem back to the case studied here.

The ND map $\mathcal{R}_\sigma$ can be approximated from the measured current and voltage data with the matrix $R_\sigma$
\begin{equation}\label{eq:NDmatrix}
R_\sigma(m,n):=\sum_{\ell=1}^L \frac{\phi_\ell^mv_\ell^n}{|e_\ell|},\qquad 1\leq m,n\leq num_{LI}
\end{equation}
where $L$ denotes the number of electrodes used, $num_{LI}$ is the number of linearly independent current patterns applied (maximum is $L-1$), and $\phi^m$, and $v^n$ denote the normalized $m$-th current pattern vector and $n$-th voltage vectors (see [\cite{Isaacson2004}; \cite{Hamilton2018_Robust}] for scaling details).  The methods described below assume the boundary conductivity $\sigma_0=1$ and that the domain has a maximum radial component of~1.  However, if this is not the case for the measured data, the ND matrix $R_\sigma$ can be scaled appropriately, as described in [\cite{Isaacson2004}], reducing the problem to the case studied here. 
\subsection{Intro to D-bar Methods for 2D EIT}\label{sec:Dbar}
While various D-bar based reconstruction algorithms for 2D EIT exist, they all have the same main structure:
\[[\textsc{Current \& Voltage Data}]\overset{1}{\longrightarrow}[\textsc{Scattering data}]\overset{2}{\longrightarrow}[\textsc{Conductivity}].\]
The scattering data is non-physical, and can be thought of as a nonlinear Fourier transform.  The D-bar methods differ in the particular formulas used to compute the scattering data and recover the conductivity.  D-bar methods come from inverse-scattering theory, an area of mathematics that brought the elegant solution to the Korteweg-de Vries (KdV) equation.  D-bar methods for EIT get their name from a $\dbar$ (D-bar) equation used to recover the conductivity $\sigma$ in Step~2 above.

Here we simulate our training data using using a variation of the `Shortcut D-bar Method' [\cite{Astala2010}] which blends the D-bar method from the Schr\"odinger equation and that of the Beltrami equation.  This is done to allow us to train the network using $L^\infty$ conductivities (Beltrami method) but still reconstruct the conductivity from the scattering data using the Schr\"odinger $\dbar_k$ equation which [\cite{Astala2010}] suggest is more robust than Step~2 of the Beltrami method.  A recent paper by \cite{Lytle2018} in fact prove that the integral equations in the Schr\"odinger formulation of the D-bar method hold for $L^\infty$ conductivities which are one near $\bndry$.

\subsubsection{Algorithm for Simulating the Training Data}\label{sec:SimDbarAlg}
Let $\Omega$ be the unit disc.  Given a set of $N$ conductivities $\left\{\sigma_n\right\}_{n=1}^N$ in $L^\infty(\Omega)$, for each $\sigma_n$ compute the associated low-pass D-bar reconstruction $\sigma_n^{\mbox{\tiny \bf DB}}$ as follows: 1) Generate the Beltrami scattering data $\tau(k)$ for $|k|\leq R$ for some chosen radius $R>0$, and 2) Solve the Schr\"odinger $\dbark$ equation using the Beltrami scattering data for $|k|\leq r$ where $r\leq R$.

\vspace{1em}
\begin{enumerate}
\item[{\underline{\it Step 1:}}] Generate the Beltrami scattering data $\tau_n(k)$ for $\sigma_n(z)$ for $k\in\C$, $|k|\leq R$ as in [\cite{Astala2010}]
 \begin{equation}\label{eq:BelScat}
 \overline{\tau_n(k)}:=\frac{1}{2\pi}\int_{\R^2}\dbarz\left[M_{+\mu_n}(z,k)-M_{-\mu_n}(z,k)\right]dz_1dz_2
\end{equation}
where $M_{\pm\mu_n}(z,k)=e^{-ikz}f_{\pm\mu_n}(z,k)$ are solutions to the Beltrami equation
\begin{equation}\label{eq:BelEq} 
\dbarz f_{\pm\mu_n}(z,k) = \pm\mu_n(z) \overline{\dez f_{\pm\mu_n}(z,k)}
\end{equation}
 satisfying $M_{\pm\mu_n}(z,k)= 1+ \mathcal{O}\left(\frac{1}{|z|}\right)$ for large $|z|$ and $\mu_n(z)=\frac{1-\sigma_n(z)}{1+\sigma_n(z)}$ denotes the corresponding Beltrami coefficient.

\vspace{1em}
\item[{\underline{\it Step 2:}}]  Relate the Beltrami and Schr\"odinger scattering data via $\T_n(k)= -4\pi i \overline{k}\tau_n(k)$, setting $\T_n(k)=0$ for all $|k|>R$.  Recover the low-pass D-bar reconstruction $\sigma_n^{\mbox{\tiny \bf DB}}=\left[m_n(z,0)\right]^2$ by solving the Schr\"odinger  $\dbark$ equation [\cite{Knudsen2009}]
\begin{equation}\label{eq:dbark}
\dbark m_n(z,k) = \frac{1}{4\pi\bar{k}}\T_n(k)e(z,-k)\overline{m_n(z,k)},
\end{equation}
for each $z\in[-1,1]^2$, where $e(z,k):=\exp\{i(kz+\bar{k}\bar{z})\}$ is a unitary multiplier, using the integral form
\begin{equation}\label{eq:dbark_intform}
m_n(z,\kappa) = 1+ \frac{1}{4\pi^2}\int_{\C}\frac{\T_n(k)e(z,-k)}{(\kappa-k)\bar{k}}\overline{m_n(z,k)}\;d\kappa_1d\kappa_2,
\end{equation}
\end{enumerate}
and the computational method outlined in [\cite{Mueller2002}].

Note that no electrode or boundary information is used in the training data as $\mu_n(z)=0$ near $\bndry$.  The choice of $\Omega=\D$ does not include boundary specific information since in the reconstruction step from experimental data, we will scale the ND map by the maximum radial component of the experimental domain $\Omega_{\mbox{\tiny meas}}$, shrinking the problem to exist within our studied domain $\Omega=\D$.   Additionally, note that the integral in \eqref{eq:dbark_intform} reduces to an integral over $|k|\leq R$ due to the compact support of $\T_n(k)$, and from [\cite{Nachman1996}] $\frac{\T_n(k)}{\bar{k}}=0$ for $k=0$.

\subsubsection{Recovery of Conductivity from Experimental Data}\label{sec:DbarExper}
Recover the D-bar reconstruction $\sigma^{\mbox{\tiny DB}}$ from the measured current and voltage data via a modification to the  Schr\"odinger $\T$ `exp' method as follows.  

\vspace{1em}
\begin{enumerate}
\item[{\underline{\it Step 1:}}] Compute the modified Schr\"odinger `exp' scattering data
\begin{eqnarray}
\texp(k) &=& \int_{\partial\Omega_1} e^{i\bar{k}\bar{z}}\left(\Lambda_{\sigma}-\Lambda_1\right){e^{ikz}}ds(z)\nonumber\\
&=& \int_{\partial\Omega_1}e^{i\bar{k}\bar{z}}\left[\Lambda_{\sigma}\left({e^{ikz}}\right)-ik\nu{e^{ikz}}\right]ds(z),\label{eq:texp}
\end{eqnarray}
for $k\in\C\setminus 0$,  $|k|\leq R_{\mbox{\tiny meas}}$ for some chosen radius $0<R_{\mbox{\tiny meas}}\leq R$.  

\vspace{1em}
\item[{\underline{\it Step 2:}}]  Recover the D-bar conductivity reconstruction $\sigma^{\mbox{\tiny DB}}=\left(m^{\mbox{\tiny exp}}(z,0)\right)^2$ using \eqref{eq:dbark} with $\texp$ in place of $\T_n$, setting $\frac{\texp(k)}{\bar{k}}=0$ for $k=0$.
\end{enumerate}

\vspace{1em}
The second line \eqref{eq:texp} comes from computing $\Lambda_1e^{ikz}= 1\nabla\left(e^{ikz}\right)\cdot\nu = ik\nu{e^{ikz}}$ which uses a continuum approximation for the DN map $\Lambda_1$ where $\nu=\nu(z)$ is the unit outward facing normal to the scaled boundary $\bndry_1$ which has maximal radial component 1.  The DN matrix approximation to $\Lambda_\sigma$ is computed from $L_\sigma=\left(R_\sigma\right)^{-1}$ via \eqref{eq:NDmatrix}.  The DN map is also scaled by the radius of the smallest circle containing the imaged domain $\Omega_{\mbox{\tiny meas}}$, and $\sigma_0$ the conductivity near the boundary $\bndry_{\mbox{\tiny meas}}$.  If $\sigma_0$ is unknown, the best constant-conductivity fit to the measured data can be used as described in [\cite{Cheney1990}].  The resulting conductivity at the end of the algorithm is then re-scaled by $\sigma_0$.  Here we compute $\nu$ numerically using a parameterization of the approximate boundary shape function (see [\cite{Hamilton2018_Robust}] for robustness studies of D-bar methods to incorrect boundary shape).  Note that we only require the measured current and voltage data, approximate boundary shape of the imaged domain $\Omega_{\mbox{\tiny meas}}$, and approximate locations of the electrodes for the D-bar reconstruction $\sigma^{\mbox{\tiny DB}}$.

\subsection{Deep Learning and U-Net}\label{sec:CNN}
In this study we follow the approach proposed in [\cite{Kang2017} and [\cite{Jin2017}] for post-processing corrupted reconstructions, in our case given by the D-bar algorithm described above in Section~\ref{sec:Dbar}. The methodology is motivated by the fact that the initial reconstruction is of convolutional type, such as the normal operator in CT, or in our case inversion of the truncated scattering transform. Consequently, we follow [\cite{Jin2017}] where the authors propose that a CNN can be used to remove artefacts and recover resolution loss present in the initial reconstruction. The network architectures used for this task are based on the well established U-Net [\cite{Ronneberger2015}], a multiscale autoencoder. 

Let us denote the network by $G_\Theta$, where $\Theta$ denotes the network parameters consisting of convolutional filters and biases. Then the learning task is an optimization problem to find the optimal set of parameters such that a loss function is minimized with respect to a training set. Specifically, in our case the training set is given by ground truth conductivities  $\sigma_n$ and D-bar reconstructions $\sigma_n^{\mbox{\tiny \bf DB}}$ for $n\in \mathcal{N}=\{1,\dots,N\}$, both given on the square $[-1,1]^2$. Then the aim is to find a network that maps from D-bar reconstruction to the correct ground truth conductivity, hence we aim to find the optimal set of parameters such that
\begin{equation}\label{eqn:l2loss}
\Theta = \arg\min_{\Theta} \sum_{n=1}^N \| G_\Theta(\sigma_n^{\mbox{\tiny \bf DB}}) - \sigma_n \|_2^2.
\end{equation}
The optimization is typically performed in subsets (batches) of training pairs $\{\sigma_n,\sigma_n^{\mbox{\tiny \bf DB}} \}_{\mathcal{I}\subset\mathcal{N}}$, rather than the whole training set.

The chosen network architectures differ slightly depending on which task, i) or ii), of the Section~\ref{sec:intro} is considered. For scenario i.) the thoracic imaging task, we employ the same network architecture as described in [\cite{Hamilton2018_DeepDbar}] as it has been shown to be specifically suited to reproduce structures in a known constrained environment. For task ii.) with minimal \emph{a priori} knowledge, an assessment of network architectures was performed and we found that adding a residual connection as in [\cite{Jin2017}] increased robustness in recovering more general shapes that were not present in the training set. In both cases we kept the filter size of the convolutional kernels as $5\times5$ and used 4 max-pool layers, as the original U-Net architecture suggests. Networks are implemented with TensorFlow in Python.

\subsection{Evaluation of the Method}\label{sec:Evaluation}
To evaluate the effectiveness of our proposed Beltrami-net method we tested it on experimental data from two different EIT machines, namely, {\sc ACT4} from Rensselaer Polytechnic Institute (RPI) [\cite{Liu2005}] and {\sc KIT4} from the University of Eastern Finland (UEF) [\cite{Kourunen2008}].  We evaluate reconstruction quality using {\it Structural SIMilarity} Indices (SSIMs) and relative $\ell_1$ and $\ell_2$ image errors.  The ground truth inclusion boundaries were extracted from photographs of the experiments.

As an additional comparison, we include 2D a-EIT reconstructions for the KIT4 data 
using a total variation (TV) regularized least squares (LS) approach.  The discretized version of the problem is
\begin{equation}
  \hat\sigma= {\rm arg} \min_{\sigma >0}\{\Vert V-U(\sigma)\Vert^2
+ \alpha TV (\sigma)\},
\label{eq.tvref}
\end{equation}
where $\sigma \in \mathbb{R}^N$ is a piece-wise constant representation of the conductivity in a set of $N$ pixels covering
the domain $\Omega$,
$U(\sigma)$ is the finite element method (FEM) based forward solver of the complete electrode model [\cite{Somersalo1992}],
$\alpha$ is the regularization parameter and
$TV(\sigma)$ is the isotropic TV functional [\cite{Rudin1992}]
\begin{equation}
\label{eq.TV.isotropic}
TV(\sigma) = 
\sum_{k=1}^N 
\sqrt{(\mathbf{D}_x\sigma)_i^2 + (\mathbf{D}_y\sigma)_i^2 + \beta},
\end{equation}
where $\mathbf{D}_x$ and $\mathbf{D}_y$ are finite dimensional approximations for the partial derivatives.
The minimization problem (\ref{eq.tvref}) is solved by using a Gauss-Newton optimization method equipped with a line search algorithm. The line search is implemented using bounded minimization such that the non-negativity $\sigma > 0$ is enforced. A detailed exposition of the method
(\ref{eq.tvref}) can be found in [\cite{Gonzalez2017}].
The regularization parameter $\alpha$ was tuned manually for the best visual quality of the reconstruction.

\subsubsection{Experimental Data}\label{sec:data}
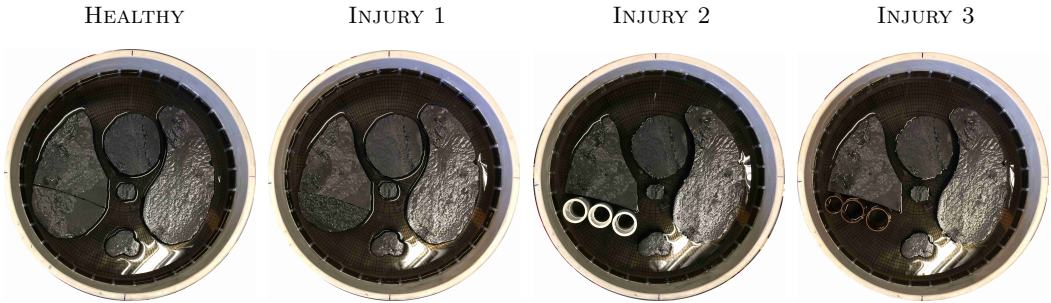
\begin{figure}[h!]
\begin{picture}(390,115)

\put(0,0){\includegraphics[width=95pt]{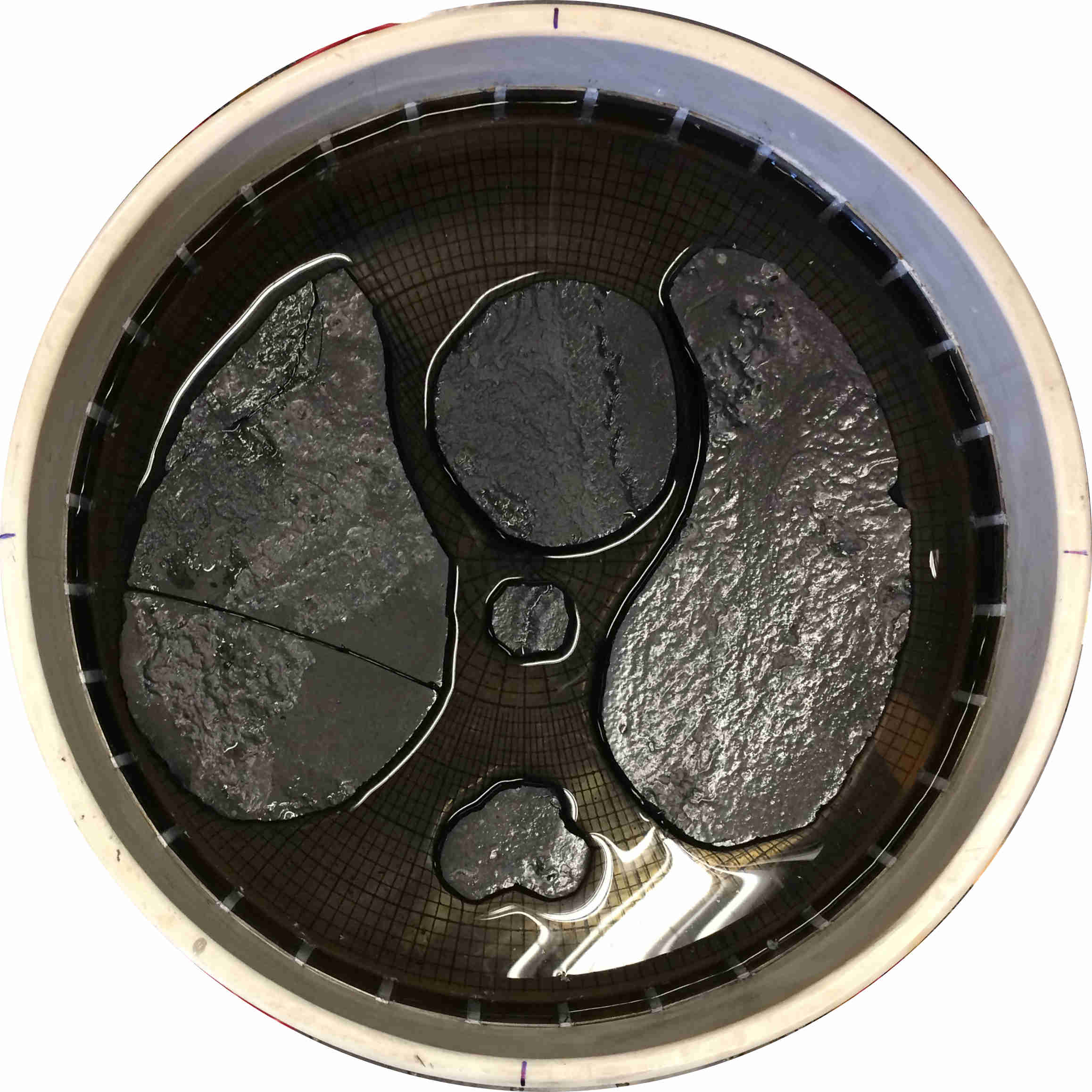}}
\put(100,0){\includegraphics[width=95pt]{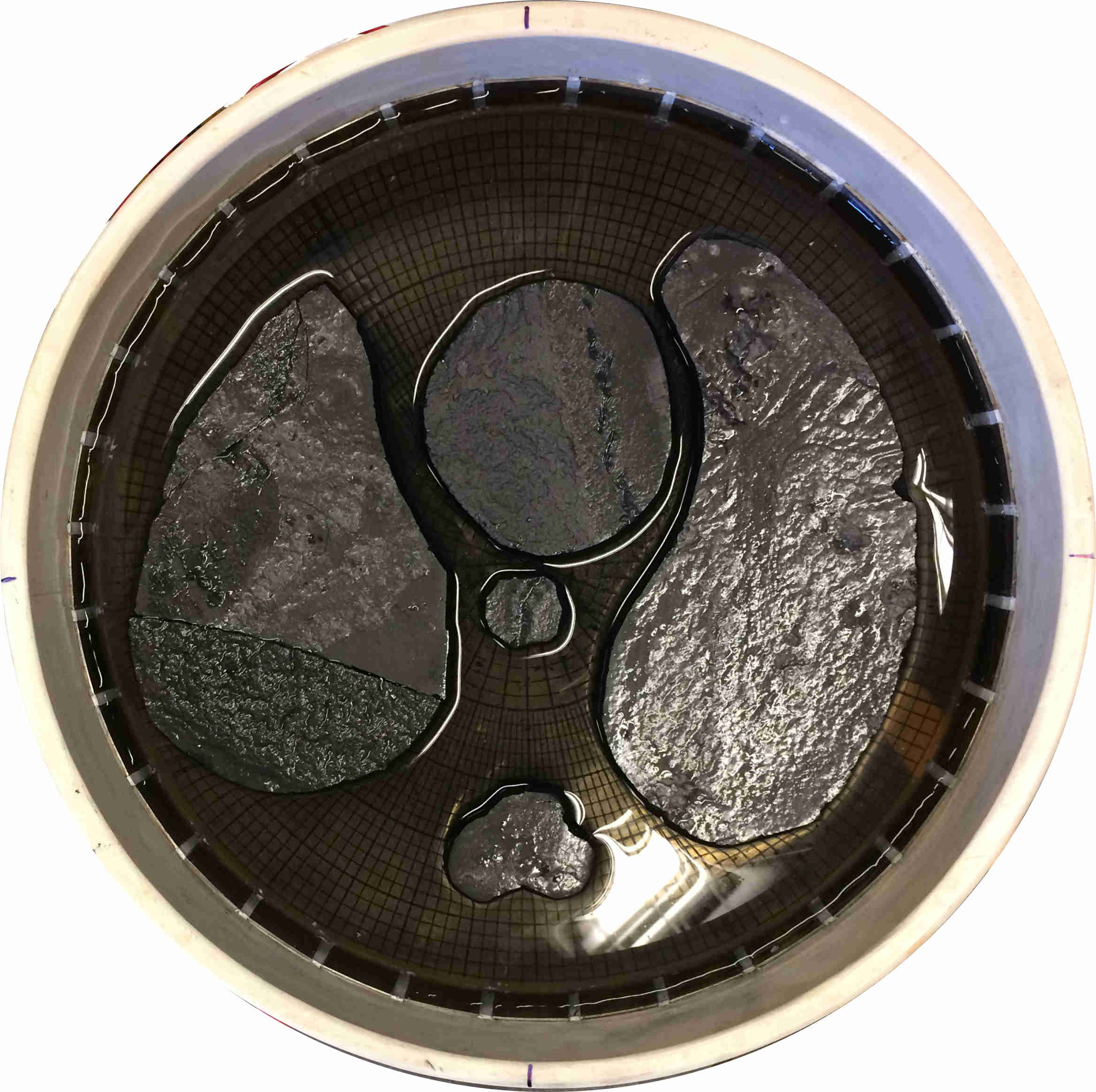}}
\put(200,0){\includegraphics[width=95pt]{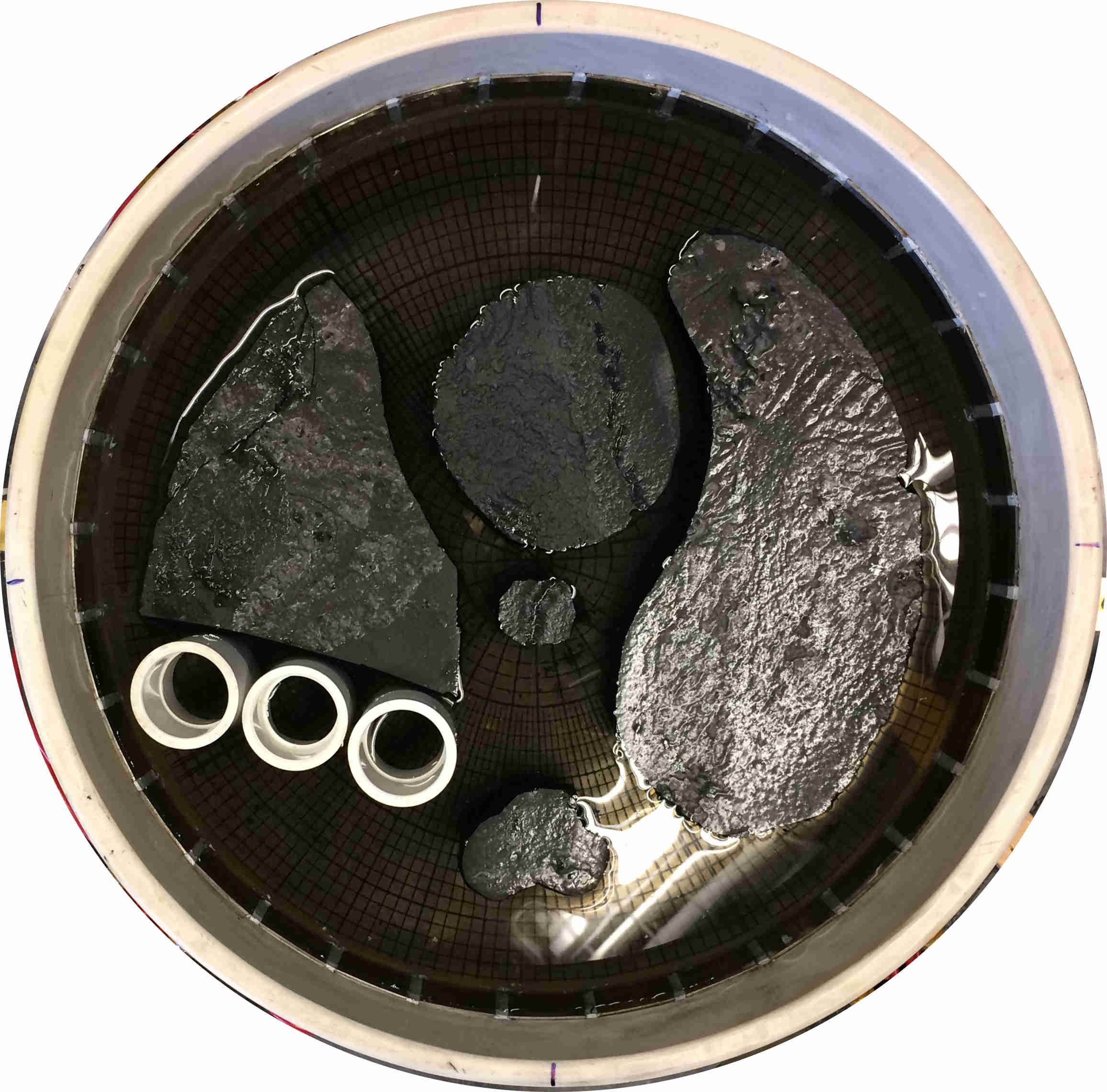}}
\put(300,0){\includegraphics[width=95pt]{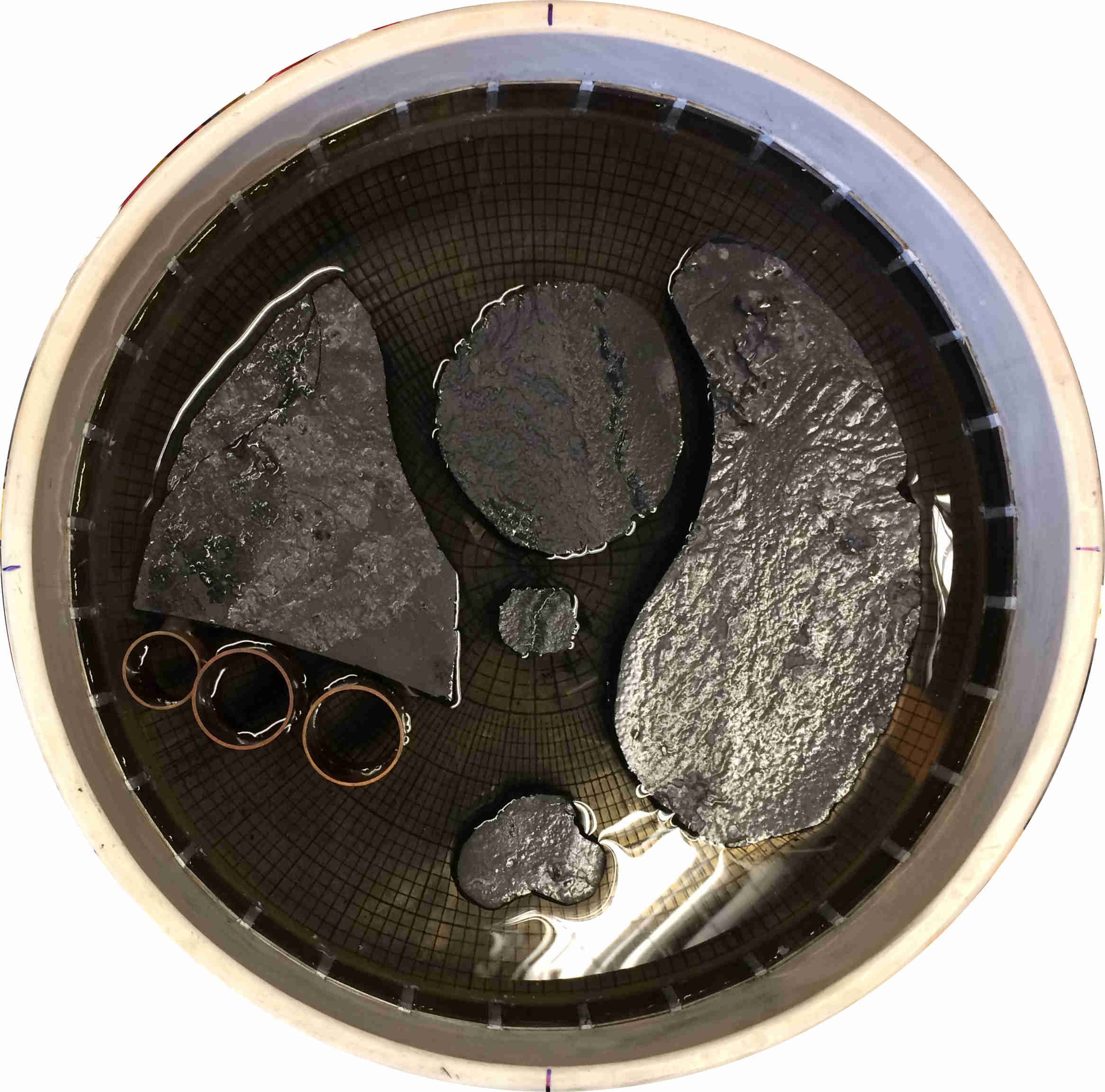}}

\put(30,105){\sc \small {Healthy}}

\put(130,105){\sc \small {Injury 1}}

\put(230,105){\sc  \small {Injury 2}}

\put(330,105){\sc \small {Injury 3}}
\end{picture}
\caption{The experimental setups for the ACT4 data collection.  Four scenarios were tested beginning with a `Healthy' setup: conductive heart and aorta, resistive lungs and spine.  In `Injury 1', the bottom portion of the right (DICOM orientation) lung was removed and replaced with a conductive agar target matching the conductivity of the heart/aorta.  In `Injury 2', the removed portion of the right lung was replaced with three plastic pipes and for `Injury 3' the removed portion is replaced with three copper pipes.}
\label{fig:RPI_phantoms}
\end{figure}

Archival ACT4 data, taken on a circular tank of radius 15cm with 32 electrodes (width 2.5cm), was used.  Agar targets with added graphite were placed in a saline bath (0.3 S/m) filled to a height of 2.25cm.  Conductive and resistive targets were used to simulate the heart and aorta, as well as the lung and spine, respectively.  See Figure~\ref{fig:RPI_phantoms} for the experimental setups.  Table~\ref{table:ACT4_setup} displays the measured conductivities of the targets, using test-cells, computed via Impedimed's SFB-7 bioimpedance meter\footnote{\url{https://www.impedimed.com/products/sfb7-for-body-composition/}}.  Trigonometric voltage patterns, with maximum amplitude 0.5V, were applied at a frequency of 3kHz and the resulting currents measured.  For consistency with previous studies, a change of basis was performed on the measured current and voltage data to synthesize the data that would have occurred if current had been applied instead of voltage (see [\cite{Hamilton2018_DeepDbar}]).  The ND and DN matrices were then computed as described in Section~\ref{sec:methods}, equation \eqref{eq:NDmatrix}.

\begin{table}[h] 
\scriptsize
  \caption{Conductivity Values for ACT4 targets at 3.3kHz}
    \begin{tabular}{l|c|c}
    \hline
&{\sc Measured Values} & {\sc Simulated Values}  \\
& (S/m) & Ranges (S/m) \\
    \hline
    \hline
    {\sc Heart/Aorta} &0.67781 & [0.5, 0.8]\\

    {\sc Lungs/Spine} & 0.056714  & [0.01, 0.2]\\
    
    {\sc Saline Background} & 0.3  & [0.29, 0.31]\\
    
    {\sc Injury 1: Agar/Graphite }&0.67781 & [0.01, 1.5]\\
            
            {\sc Injury 2: Plastic Tubes} &0  & [0.01, 1.5]\\
            {\sc Injury 3: Copper Tubes} &infinite  & [0.01, 1.5]\\
    \hline
    \hline
    \end{tabular}%
  \label{table:ACT4_setup}%
\end{table}%

\begin{figure}
\begin{picture}(410,115)
\put(95,0){\includegraphics[height=85pt]{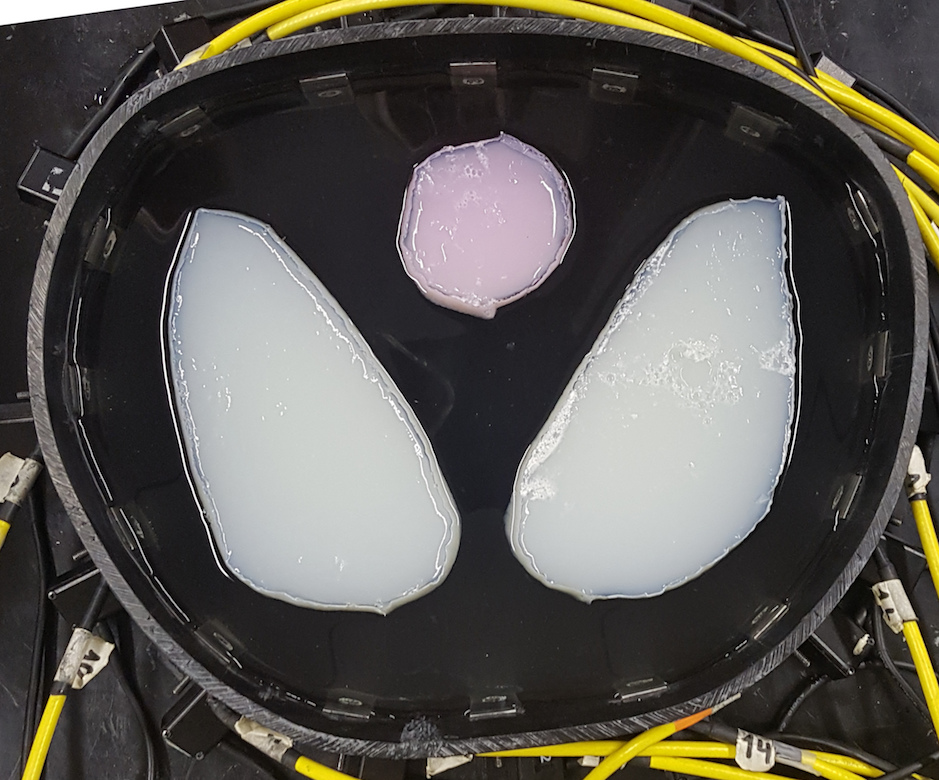}}
\put(205,0){\includegraphics[height=85pt]{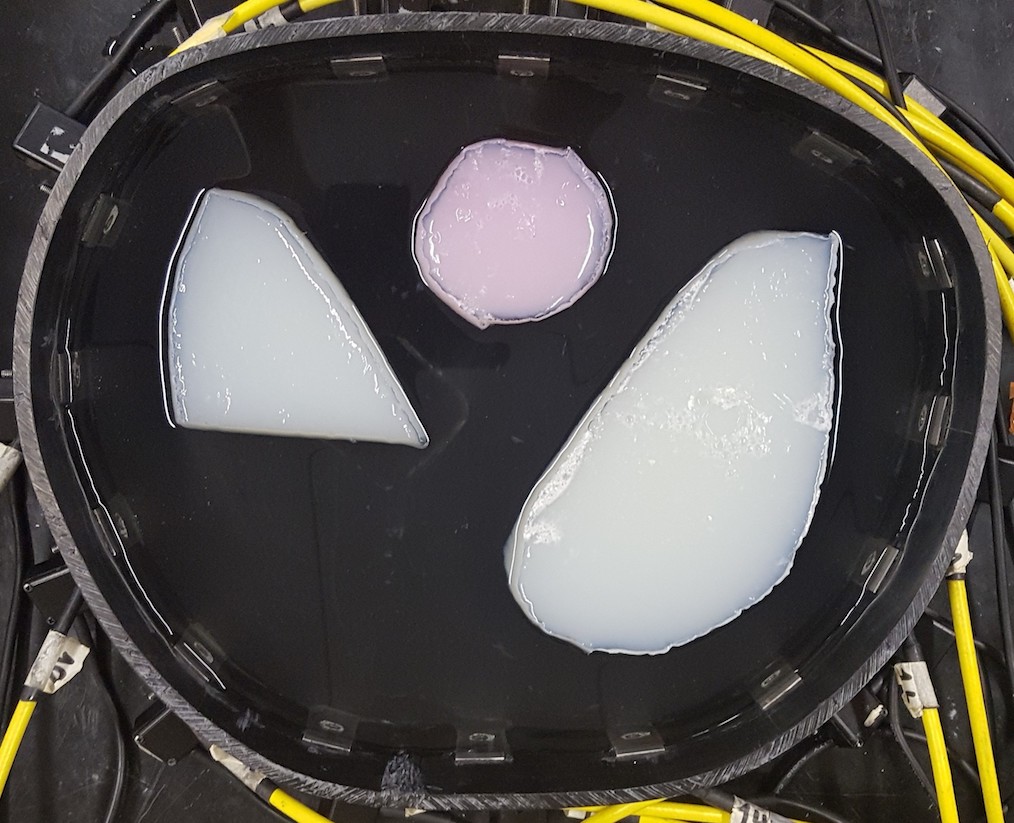}}
\put(315,0){\includegraphics[height=85pt]{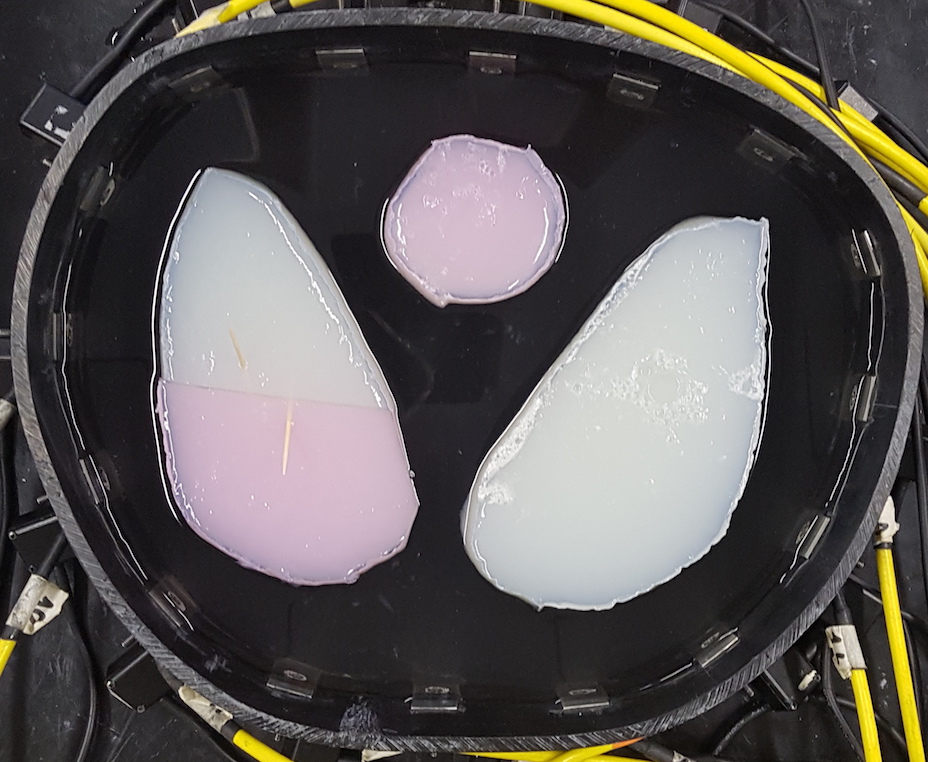}}
\put(0,0){\includegraphics[height=85pt]{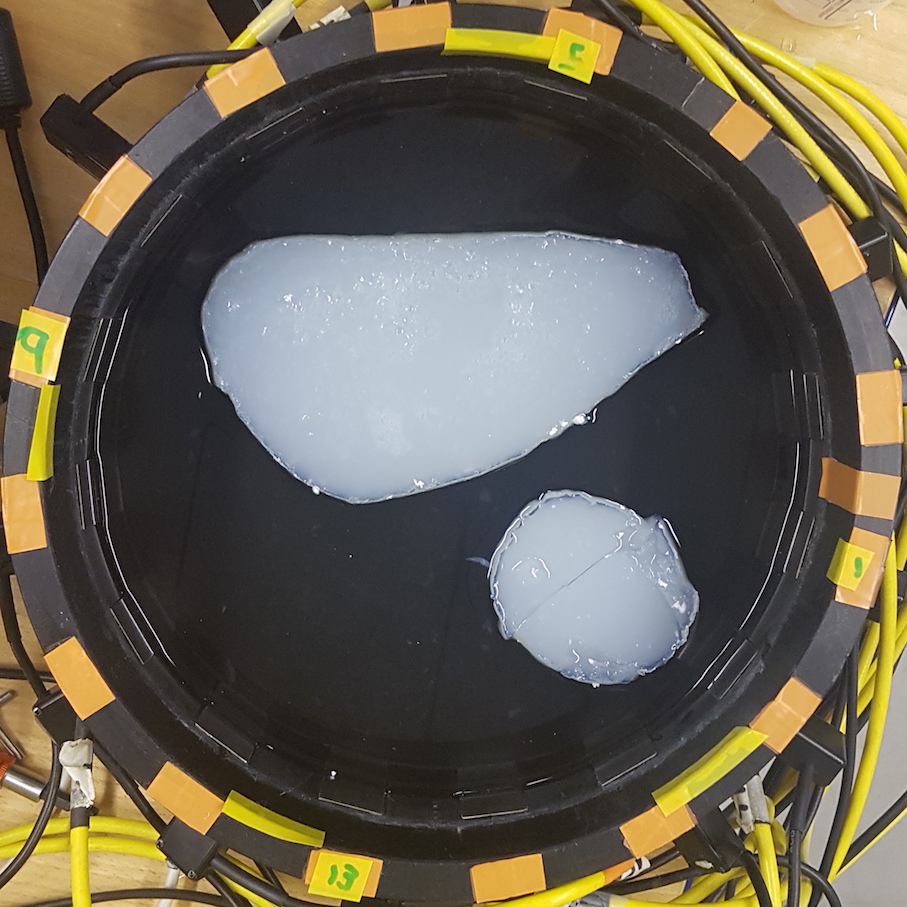}}

\put(30,95){\small {\sc Circle}}
\put(115,95){\small {\sc Chest-Healthy}}
\put(235,95){\small {\sc Chest-Cut}}
\put(340,95){\small {\sc Chest-Split}}



\end{picture}
\caption{Experimental setups for the KIT4 data on three different experimental tank setups. Circle: The large object is low conductivity and small object is high conductivity.  Chest: The agar targets are either high (pink) or low (white) conductivity.}
\label{fig:KIT4setups}
\end{figure}
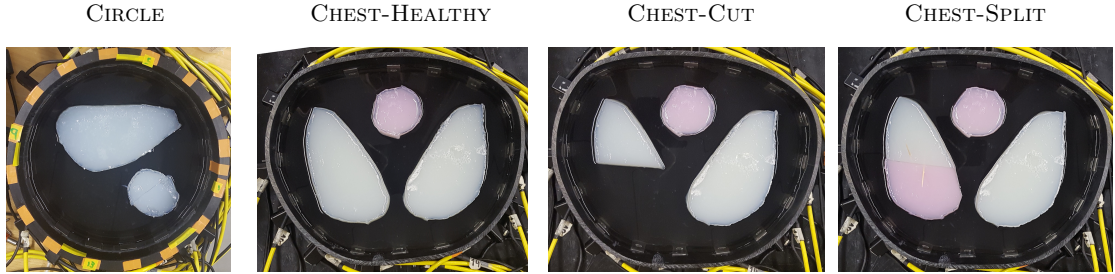

We collected KIT4 data using two different, translationally symmetric tanks 
to obtain data for two different boundary shapes, namely circle and chest-shaped, as shown in Figure~\ref{fig:KIT4setups}.  
In each tank, the number of electrodes is sixteen.
Adjacent (skip-0) current patterns were applied with current frequency at 10kHz and amplitude 3mA.  Conductive 
and resistive 
agar targets were used across all the KIT4 experiments. %
The circular tank has a radius of 14cm with 16 electrodes of width 2.5cm.  Agar targets of conductivity 67~mS/m (large object on the top) and 305~mS/m (smaller, nearly circular object on the bottom right) were placed in a saline bath of conductivity 135~mS/m filled to a height of 45mm.  
The chest shaped tank has a perimeter of 1.02m with 16 electrodes of width 2cm attached.  The locations of the electrodes are not exactly equidistant from one another but can be seen from the photographs (see Figure~\ref{fig:KIT4setups}).  Agar targets consisting of high conductivity 323~mS/m (targets with pink ink) and low conductivity 61~mS/m (white) were placed in a saline bath (conductivity 135~mS/m, height 47mm for the {\it Chest-Healthy}and {\it Chest-Cut} targets, and 44mm for the {\it Chest-Split} target in Figure~\ref{fig:KIT4setups}). The right (DICOM) lung was cut and two simulated injuries explored: 1) the bottom portion was removed completely (Fig.~\ref{fig:KIT4setups}: {\it Chest-Cut}) and 2) the bottom portion was replaced with a higher conductivity piece of agar (Fig.~\ref{fig:KIT4setups}: {\it Chest-Split}).  
%

\subsection{Training Data}\label{sec:training}
Two sets of training data were used in this study, tailored to the ACT4 and KIT4 experiments.  We introduce the notation $\widetilde{\sigma}$ to denote a conductivity that has not yet been scaled to a boundary conductivity of 1, reserving  $\sigma$ solely for conductivities with a boundary value of 1.

\subsubsection{ACT4 phantoms}
\begin{figure}[h!]
\begin{picture}(450,120)

\put(0,0){\includegraphics[height=95pt]{HLSA_saline_DeepDbar_DICOM_smaller.jpg}}
\put(110,0){\includegraphics[height=95pt]{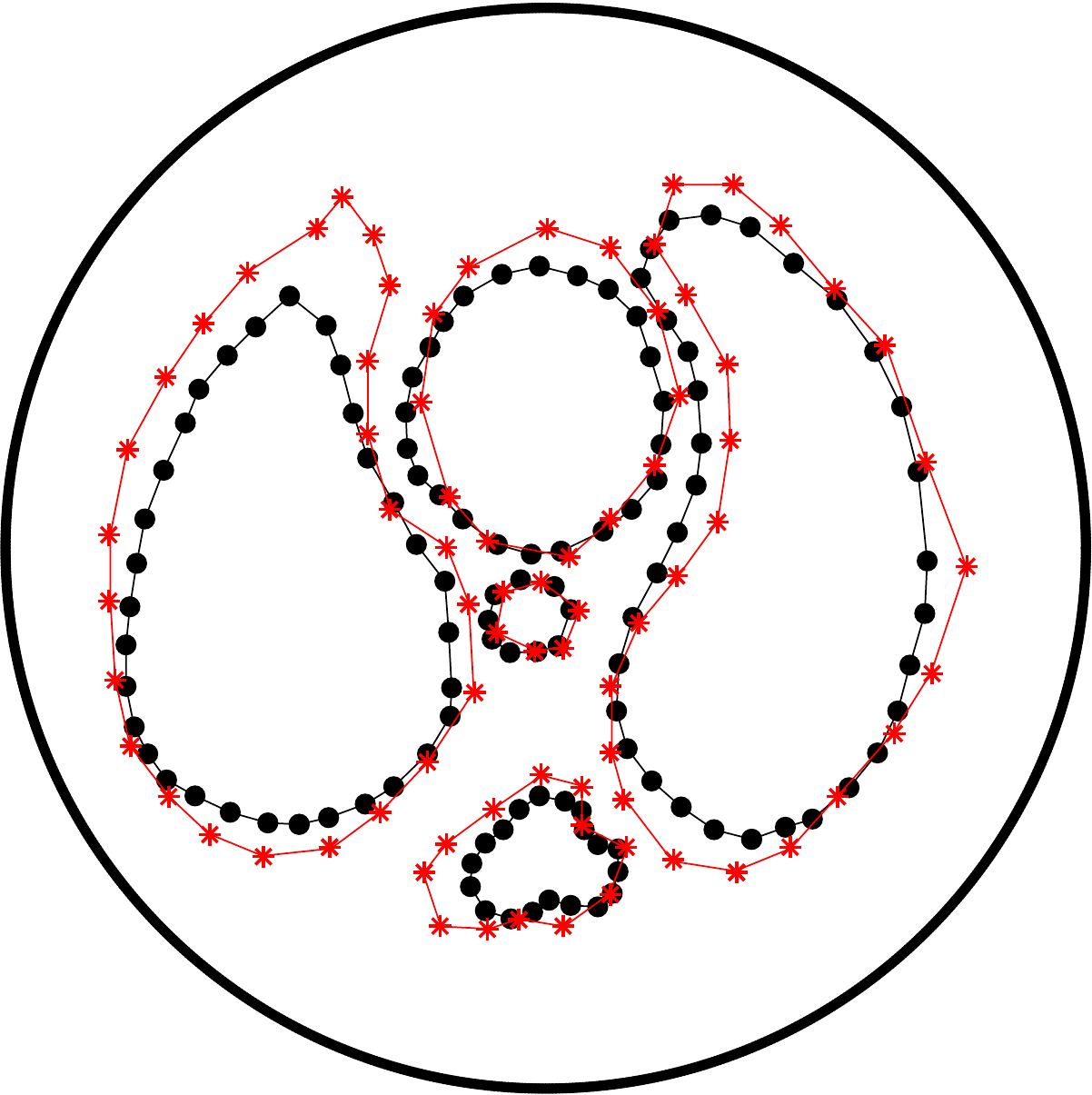}}
\put(220,0){\includegraphics[height=95pt]{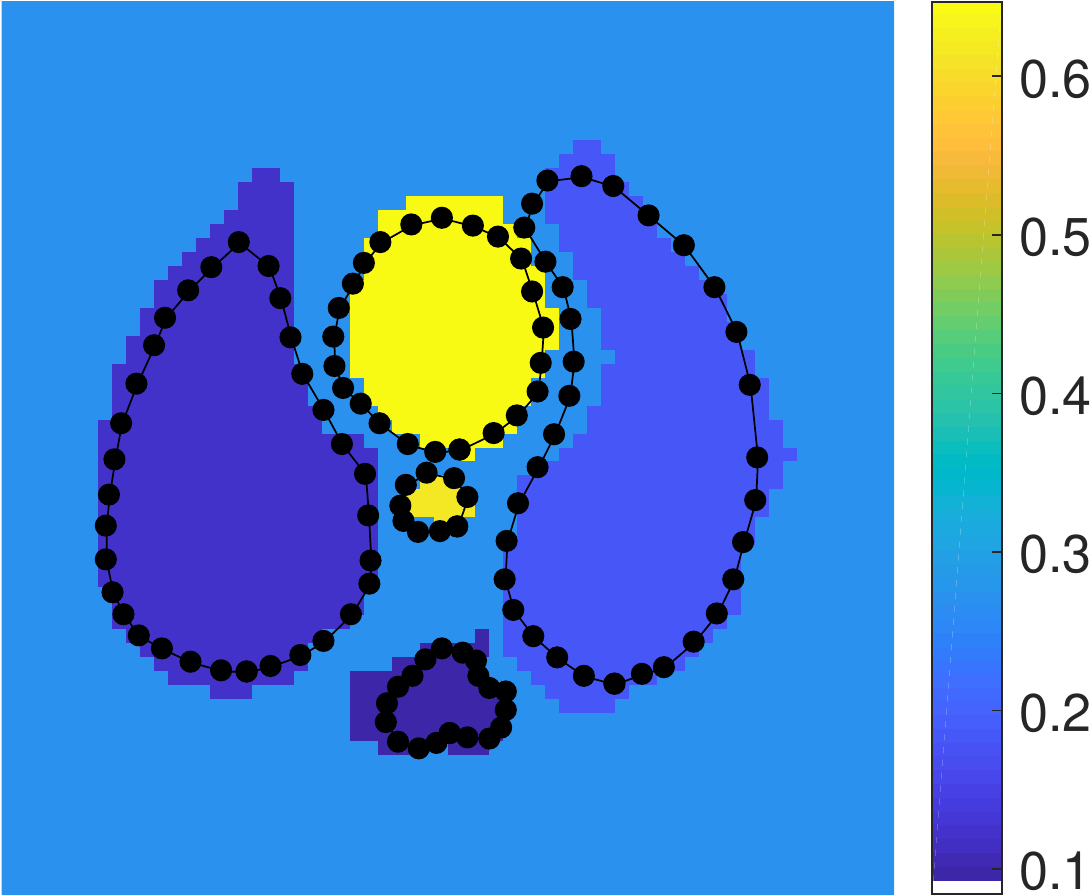}}
\put(340,0){\includegraphics[height=95pt]{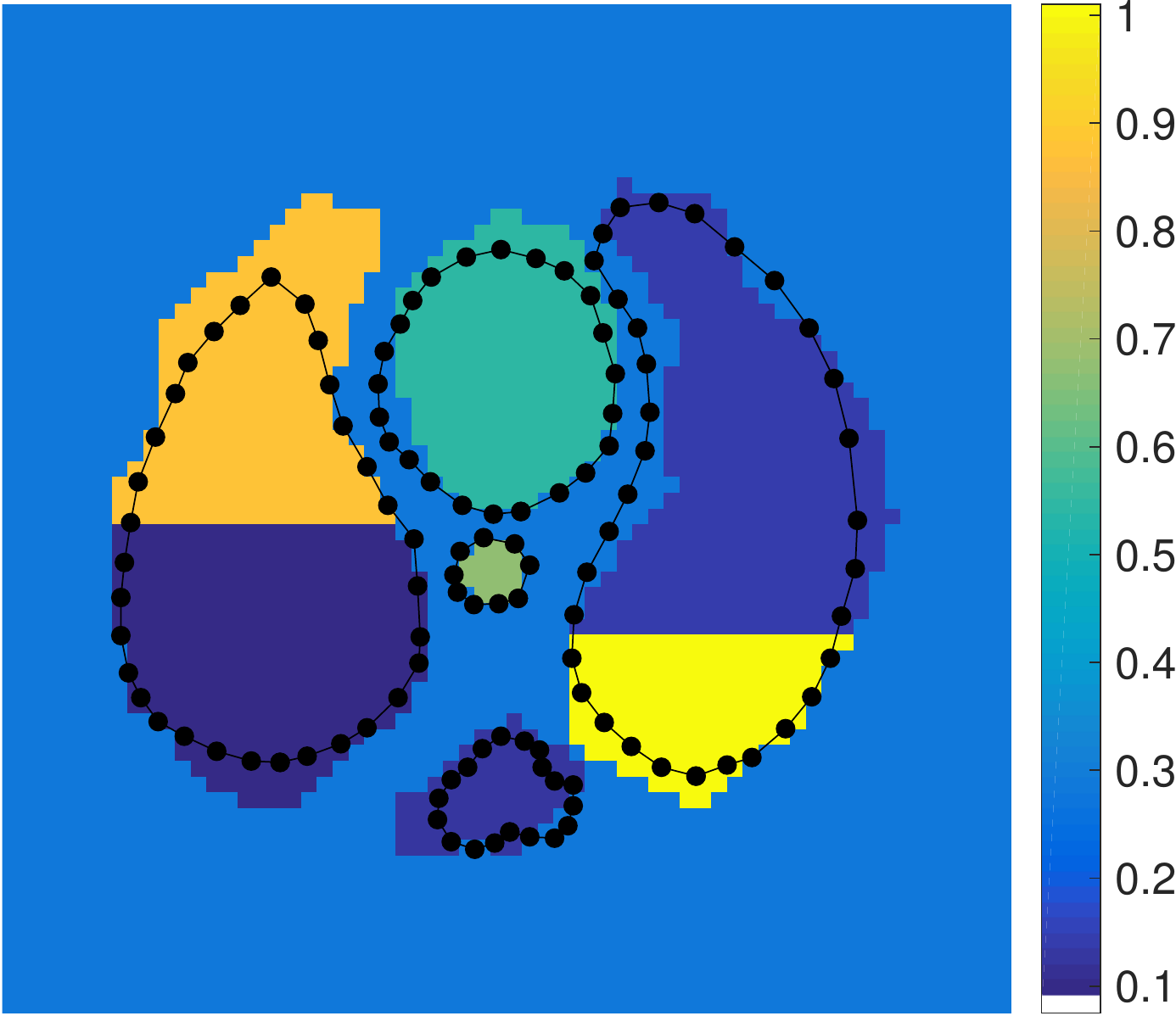}}

\put(20,110){\sc \footnotesize{{ACT4 Healthy}}}

\put(115,115){\sc \footnotesize {True \& Approximate}}
\put(135,105){\sc \footnotesize {Boundaries}}

\put(235,110){\sc \footnotesize{{Sample Healthy}}}

\put(355,110){\sc \footnotesize{{Sample Injured}}}

\end{picture}
\caption{Samples of the simulated conductivities used to generate the ACT4 training data corresponding to the experiments shown in Figure~\ref{fig:RPI_phantoms}.  Starting with a healthy setup (left), the `true organ boundaries' (shown in black dots) were extracted from the photograph along with an `approximate organ boundaries' (red stars) which are displayed in the second image.  Noise was added to these approximate boundary points to generate the organ boundaries used in the simulated conductivities.  Samples of such conductivities are shown in the third and fourth images with the true organ boundaries outlined in black dots.}
\label{fig:ACT4demoSims_noisyOrgans}
\end{figure}

Candidate phantoms $\widetilde{\sigma}_n$ for the ACT4 training were formed by extracting the approximate boundaries of the inclusions from the `Healthy' setup shown in Figure~\ref{fig:ACT4demoSims_noisyOrgans} (first).  The approximate boundaries are show in in red \tredbf{$\ast$} and the true boundaries are shown in black dots (Figure~\ref{fig:ACT4demoSims_noisyOrgans}, second).  Phantoms $\widetilde{\sigma}_n$ were generated as follows.  
\begin{itemize}
\item {\it Determine which objects are included.}  Random numbers were generated from the uniform distribution on $[0,1]$ to determine whether each inclusion (left lung: 90\%, right lung: 90\%, spine: 100\%, heart: 95\%, aorta: 95\%) was included in $\widetilde{\sigma}_n$.  
\item {\it Determine the conductivities of each target in $\widetilde{\sigma}_n$.}  The conductivities were assigned by drawing random numbers from uniform distributions using the respective conductivity windows outlined in Table~\ref{table:ACT4_setup}.
\item {\it Determine the locations of each target in $\widetilde{\sigma}_n$.}  The coordinates of the each inclusion were created by adding noise, using the {\tt awgn} command in {\sc Matlab}, to the `approximate' coordinates (red stars) of the corresponding inclusion, see Figure~\ref{fig:ACT4demoSims_noisyOrgans}.  
\end{itemize}
As the ACT4 experiments contained `injuries' to the right (DICOM) lung, simple injuries were simulated in the training data as follows.  For each included lung, do the following: 
\begin{itemize}
\item {\it Determine if the given lung contains an injury}. Generate a random number to determine whether or not an injury took place in the lung (50\% chance).  
\item {\it If yes, divide the lung into two regions.}.  Create a horizontal dividing line randomly by using the max and min vertical $x_2$ coordinates of the lung dividing the lung into two regions.  \item {\it Assign the injury.}  Draw a random number to determine which region (top or bottom) the `injury' took place (50-50 chance), and another random number drawn from the uniform distribution on the interval $[0.01,1.5]$ to determine the conductivity of the injured region.  
\end{itemize}
More complicated injuries were not considered here to allow for direct comparison to the previous study [\cite{Hamilton2018_DeepDbar}].  Sample phantoms $\sigma_n$ can be seen in Figure~\ref{fig:ACT4demoSims_noisyOrgans}, third and fourth images.

\subsubsection{KIT4 phantoms}
Conductivity phantoms $\widetilde{\sigma}_n$ for the KIT4 training data were more general as the sizes and locations of the targets in the experiments varied greatly.  Phantoms consisted of one to three ellipses of varying size (semi-major and minor axes chosen from the uniform distribution on $[0.2,\;0.35]$), location $\rho e^{i\theta}$ for $\rho\in[0,0.6]$ and $\theta\in[0,2\pi)$, and angular orientation in $[0,2\pi)$.  The ellipses were not permitted to overlap, and were all forced to be completely contained inside a $z$-disc of radius~0.95.  The background conductivity was chosen from the uniform distribution on the interval $[0.13,\;0.145]$.  For each inclusion, a random number was drawn to determine whether the inclusion was more or less conductive than the background (50-50 chance) and conductivities randomly assigned from the corresponding uniform distributions $[0.29,0.34]$ and $[0.05,0.075]$.  The chance of a target being split into two pieces was 1 in 3.  If split, no region could be smaller than 1/4 the size of the whole inclusion, and the split could be along any dividing line (horizontal, diagonal, vertical).  Divided inclusions were forced to either 1) have one part match the conductivity of the background, or 2) be split into a portion that is more conductive than the background and a portion that is less conductive than the background.  Sample simulated conductivities $\widetilde{\sigma}_n$ are shown in Figure~\ref{fig:KIT4demoSims}.

\begin{figure}[h!]
\begin{picture}(460,90)

\put(0,0){\includegraphics[height=80pt]{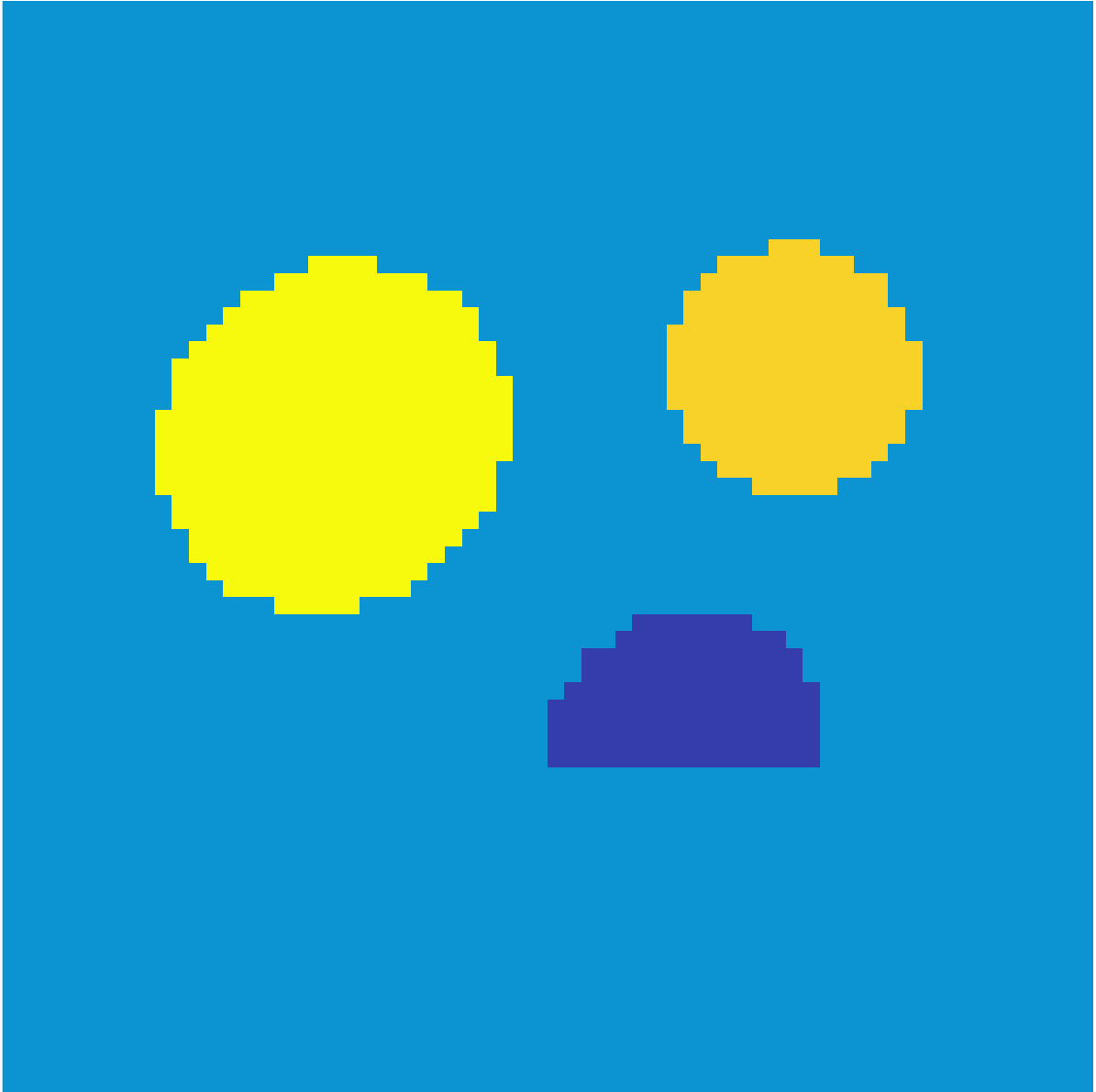}}
\put(95,0){\includegraphics[height=80pt]{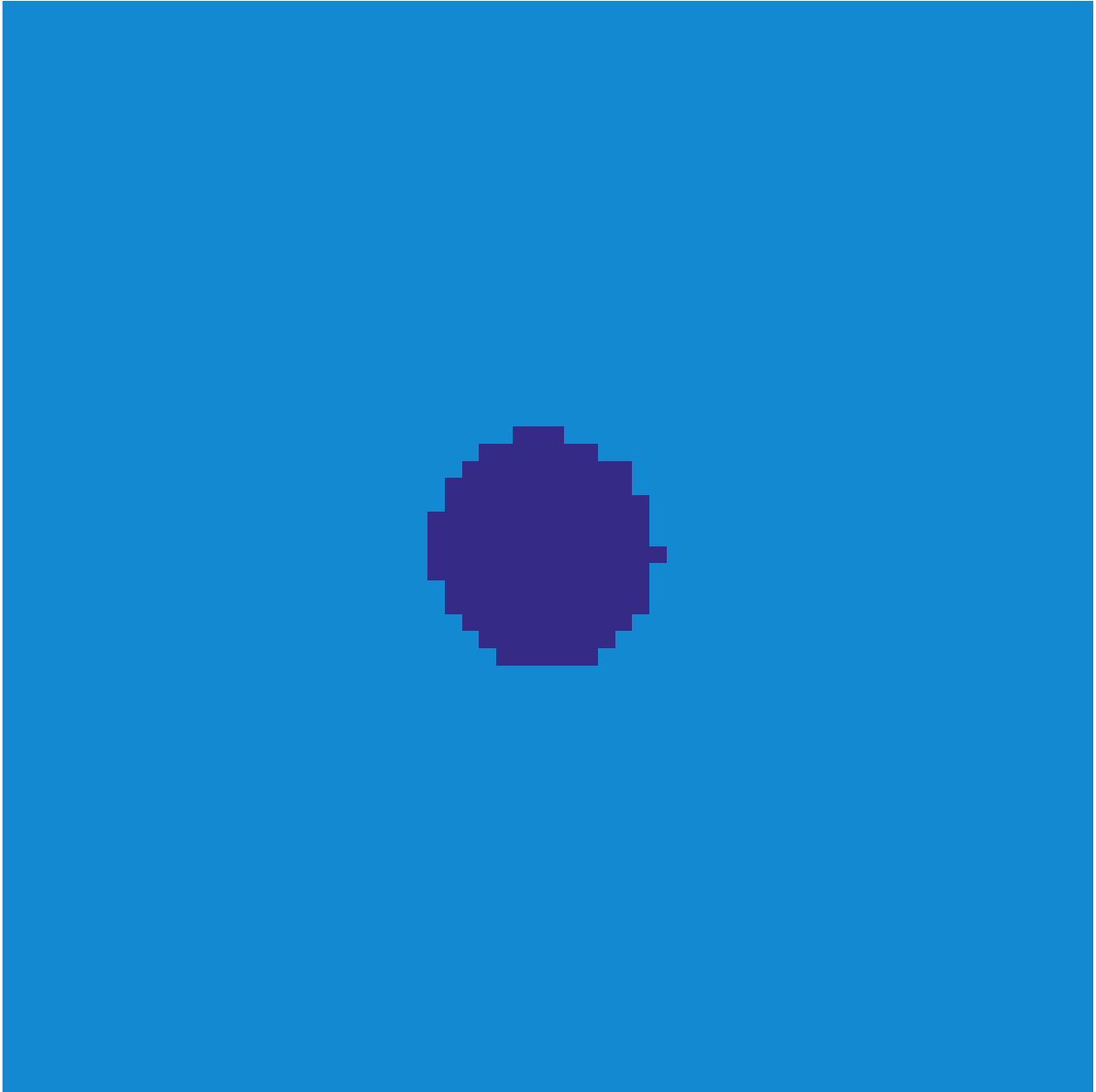}}
\put(190,0){\includegraphics[height=80pt]{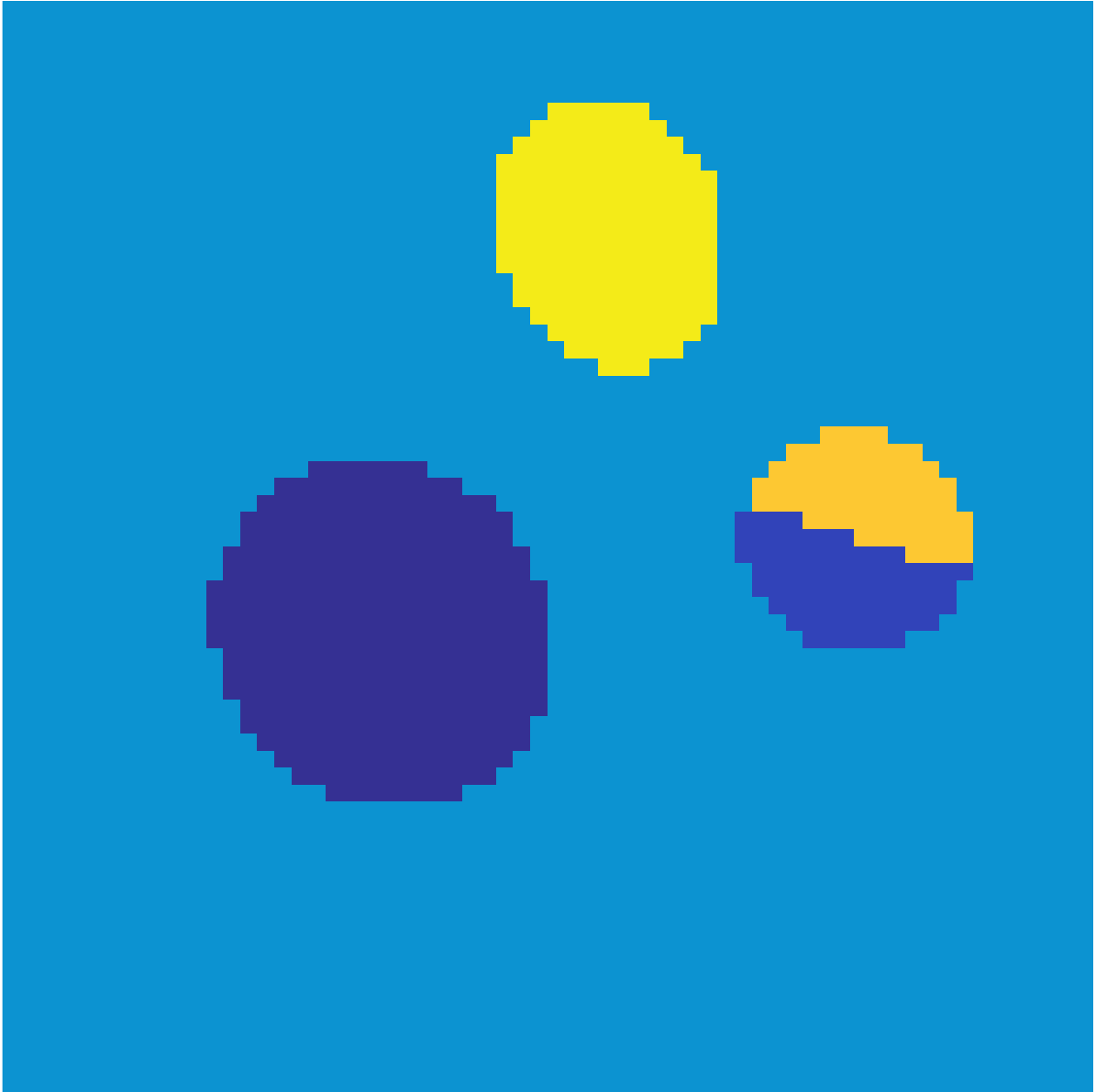}}
\put(285,0){\includegraphics[height=80pt]{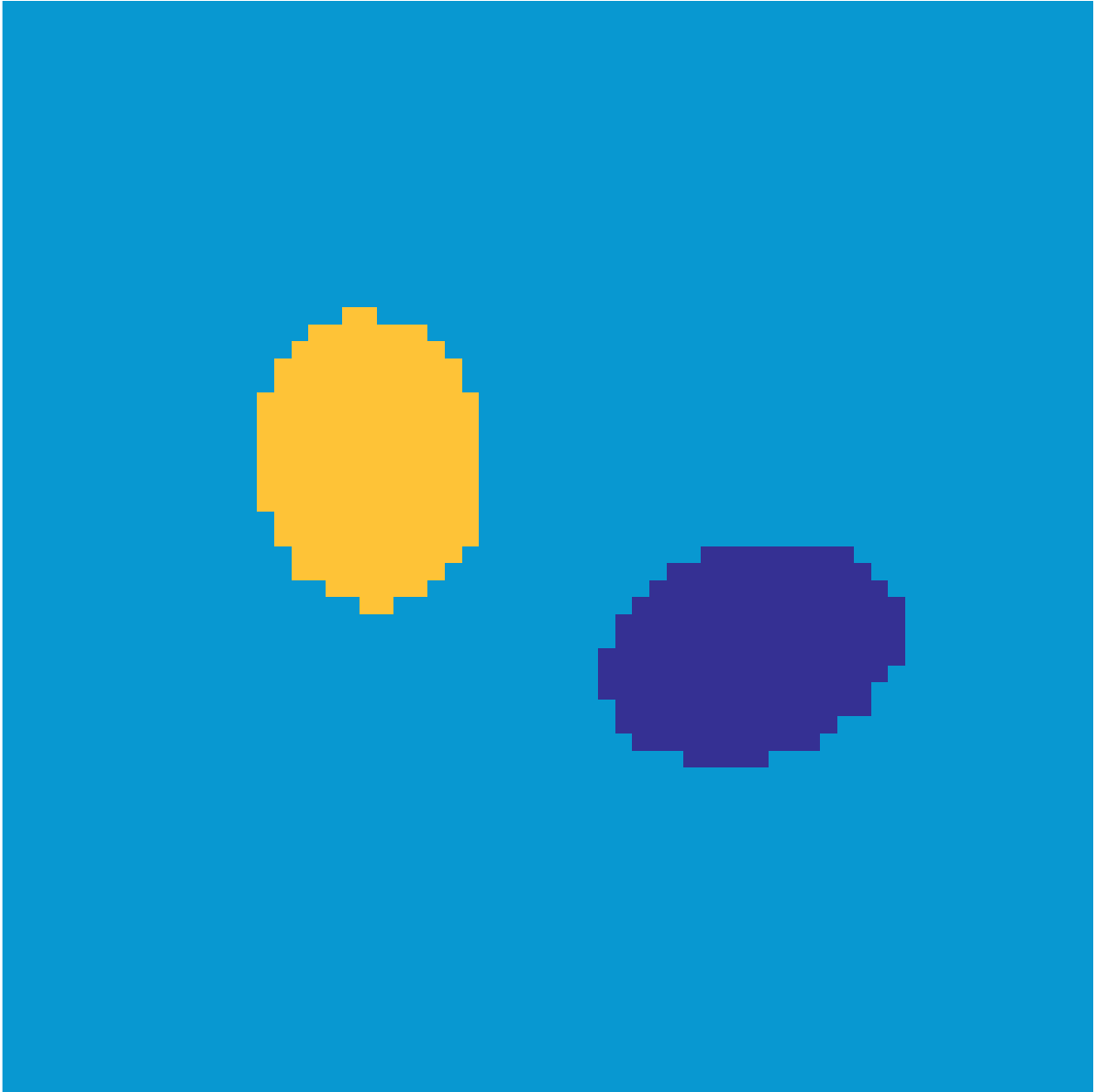}}
\put(380,0){\includegraphics[height=80pt]{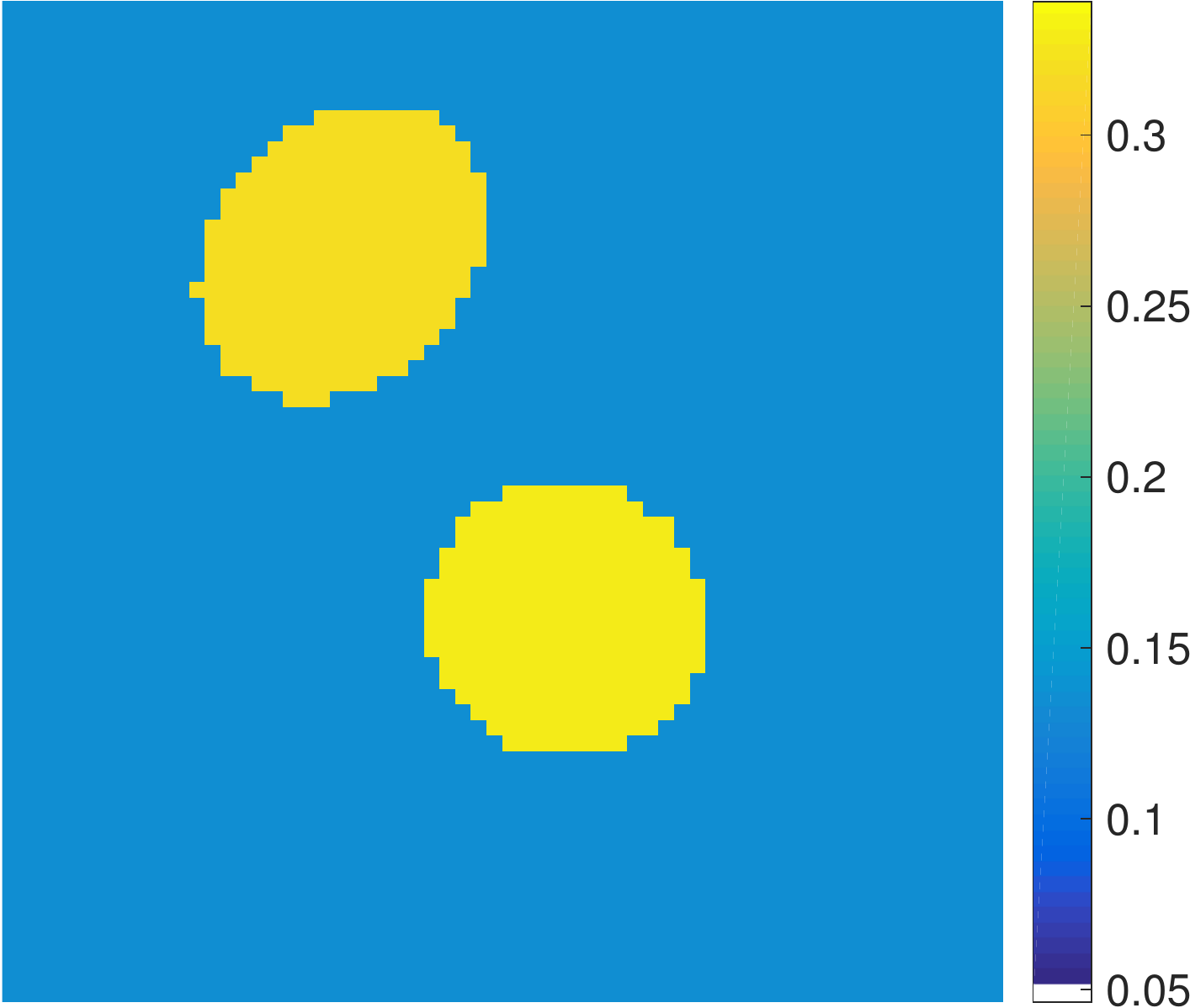}}

%
%
%
\end{picture}
\caption{Samples of the simulated conductivities used to generate the KIT4 training data corresponding to the experiments shown in Figure~\ref{fig:KIT4setups}.  One to three ellipses of varying eccentricities were randomly included with the possibility of inclusions being divided into two pieces of with no portion smaller than 1/4 of the original inclusion.}
\label{fig:KIT4Sims}
\end{figure}

\subsubsection{Producing training data}\label{sec:ProdTraining}
For each conductivity phantom $\widetilde{\sigma}_n$, the conductivity was scaled to a boundary value of~1 via $\sigma_n = \frac{1}{\sigma_{b_n}}\widetilde{\sigma}_n$ where $\sigma_{b_n}$ denotes the constant conductivity near the the boundary, here the constant background value.  If using a more complicated anatomical atlas, the value for $\sigma_{b_n}$ would be the constant conductivity for the tissue at the patient's boundary.  Then, the conductivity is extended to $[-1,1]^2$ by setting $\sigma_n=1$ for $z\in[-1,1]^2\setminus\Omega_n$.  Then, for each scaled conductivity $\sigma_n$, the Beltrami scattering data $\tau_n(k)$ \eqref{eq:BelScat} was computed for $|k|\leq R_{\mbox{\tiny ACT4}}=5$ or $|k|\leq R_{\mbox{\tiny KIT4}}= 5.5$, using a $2^5\by2^5$ uniformly spaced $k-$grid on $[-5,5]^2$ or $[-5.5,5.5]^2$, respectively, by solving \eqref{eq:BelEq} with Beltrami coefficient $\mu_n(z)=\frac{1-\sigma_n(z)}{1+\sigma_n(z)}$ as outlined in Step~1 of Section~\ref{sec:SimDbarAlg}.  Next, the blurred D-bar reconstruction $\sigDBn$ was recovered by Step~2 of Section~\ref{sec:SimDbarAlg} as follows.  First, the Beltrami $\tau_n$ was related to the Schr\"odinger $\T_n$ scattering data by $\T_n(k)= -4\pi i \overline{k}\tau_n(k)$.  Then, a random number $R_n$ was generated for the new scattering radius cutoff from the uniform distribution on $[3.5,5]$ for ACT4, or $[4, 5.5]$ for KIT4.  Then, the computed scattering data $\T_n$ was interpolated to a new $2^6\by2^6$ $k-grid$ with maximum radius $R_n$ on $[-R_n,R_n]^2$.  A non-uniform cutoff threshold was enforced by setting $\T_n(k)=0$ if $|Re(\T_n(k)|$ or $|Im(\T_n(k)|$ exceeded $thresh=24$ or $|k|>R_n$.  Then, the $\dbark$ equation was solved using the integral form \eqref{eq:dbark_intform} and the D-bar conductivity recovered as $\sigDBn(z)=\sigma_{b_n}\left(m_n(z,0)\right)^2$, rescaling by the boundary conductivity $\sigma_{b_n}$, using a $2^6\by2^6$ $z-$grid on $[-1,1]^2$ with gridsize $h_z\approx 0.0317$.

A total of $4,096$ (ACT4) and $15,360$ (KIT4) pairs $\left\{\widetilde{\sigma}_n,\sigDBn\right\}$ were created for use as training data in the U-net architectures described above in Section~\ref{sec:CNN}. Training was performed with the Adam optimizer and an initial learning rate of $10^{-4}$ to minimize the $\ell^2$-loss \eqref{eqn:l2loss} with a batch size of 16 and for a total of 200,000 iterations. Training was supervised with a simulated validation set of $\sim 5\%$ of the training set size. The long training time, in terms of iterations, was mainly necessary to obtain constant areas in the inclusions as well as background. The training procedure took roughly 3 hours for each experiment on a single Nvidia Titan XP GPU.

Then, after the successful training procedure, the effectiveness was evaluated on simulated datasets $\sigDBn$ not used in the training or validation data (Section~\ref{sec:resultsSim}) as well as experimental reconstructions for the ACT4 and KIT4 data, applied to the respective ACT4 or KIT4 network (Section~\ref{sec:resultsExp}).

\section{Results \& Discussion}\label{sec:results}
Here we present the results of the new Beltrami-Net method on experimental, as well as simulated, data from the ACT4 and KIT4 EIT systems.

\subsection{Reconstructions from Simulated Data}\label{sec:resultsSim}
We begin by visually testing the quality of the Beltrami-Net approach on simulated data. We explore test cases consistent with the training data, as well as phantoms that deviate from the procedure for creating the training set. 

Figure~\ref{fig:ACT4_Results} shows sample low-pass D-bar and Beltrami-Net reconstructions from simulated test data for the ACT4 scenario.  As it can be seen, if the injuries are consistent with the training, at most a single horizontal dividing line in the lung as in Sims 1-2, the network can almost perfectly recover the targets. If the test data deviates from this convention, Sims 3-5, it is more difficult to recover the correct location and structure, most notably for vertical divisions. Nevertheless, for two dividing lines the network is able locate the conductivity change correctly and establishes a sharp division in the reconstruction.

\begin{figure}[h!]
\begin{picture}(270,470)

\put(10,360){\includegraphics[height=80pt]{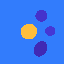}}
\put(100,360){\includegraphics[height=80pt]{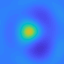}}
\put(190,360){\includegraphics[height=80pt]{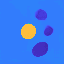}}

\put(10,270){\includegraphics[height=80pt]{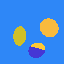}}
\put(100,270){\includegraphics[height=80pt]{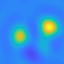}}
\put(190,270){\includegraphics[height=80pt]{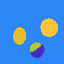}}

\put(10,180){\includegraphics[height=80pt]{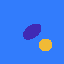}}
\put(100,180){\includegraphics[height=80pt]{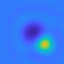}}
\put(190,180){\includegraphics[height=80pt]{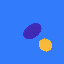}}

\put(10,90){\includegraphics[height=80pt]{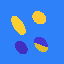}}
\put(100,90){\includegraphics[height=80pt]{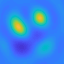}}
\put(190,90){\includegraphics[height=80pt]{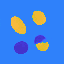}}

\put(10,0){\includegraphics[height=80pt]{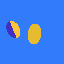}}
\put(100,0){\includegraphics[height=80pt]{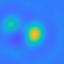}}
\put(190,0){\includegraphics[height=80pt]{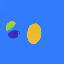}}

\put(38,455){\footnotesize \sc Truth}

\put(121,460){\footnotesize \sc Low Pass}
\put(115,450){\footnotesize \sc D-bar Image}

\put(203,460){\footnotesize \sc  Beltrami-Net}
\put(218,450){\footnotesize \sc  Image}

\put(0,30){\rotatebox{90}{{\sc \footnotesize Sim 5}}}
\put(0,120){\rotatebox{90}{{\sc \footnotesize Sim 4}}}
\put(0,210){\rotatebox{90}{{\sc \footnotesize Sim 3}}}
\put(0,300){\rotatebox{90}{{\sc \footnotesize Sim 2}}}
\put(0,390){\rotatebox{90}{{\sc \footnotesize Sim 1}}}

%
%
%
%
%
%
%
%
%

\end{picture}
\caption{Results for simulated test data with the network trained for KIT4. Note that the training data only included up to 3 inclusions. All images are on the same scale.  }  
\label{fig:KIT4demoSims}
\end{figure}
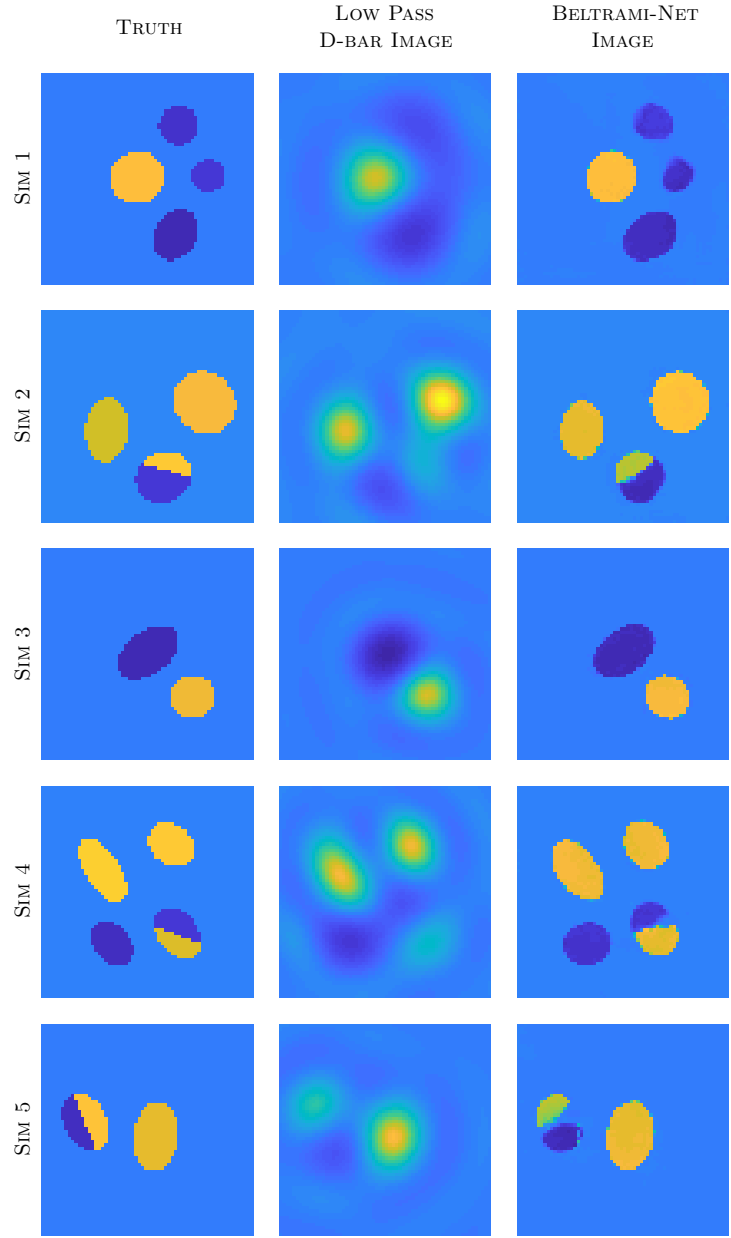

Reconstructions from simulated test data for KIT4 are shown in Figure~\ref{fig:KIT4demoSims}. Most notably, if the inclusions are isolated and do not include a cut, the network can reconstruct these very well. We note here that the training data only included up to 3 inclusions. Nevertheless, the network seems to have no difficulties to reconstruct 4 inclusions in the image. As can be seen, the cut ellipses are more difficult to reconstruct. In most cases the network manages to include a cut in the ellipse, but in a wrong orientation. In some cases, such as simulation 5, the network is not able to distinguish between a cut and two separate inclusions.

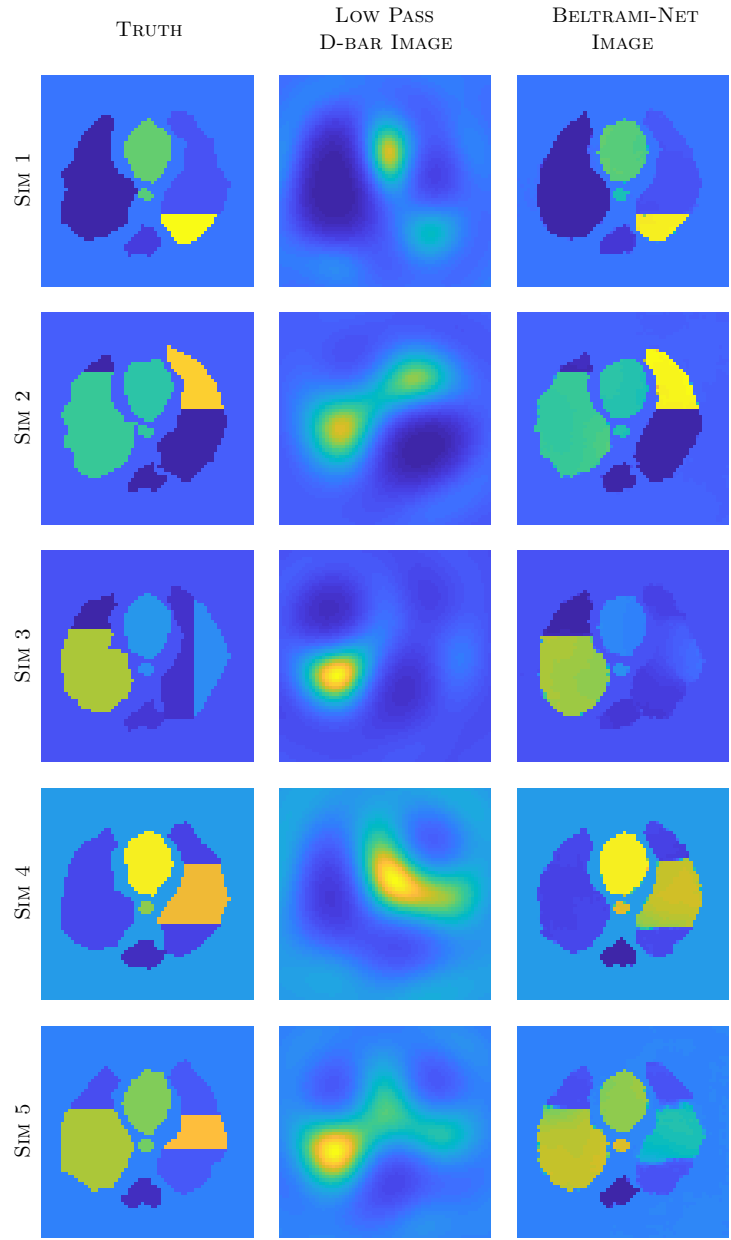
\begin{figure}[h!]
\begin{picture}(270,470)

\put(10,360){\includegraphics[height=80pt]{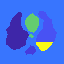}}
\put(100,360){\includegraphics[height=80pt]{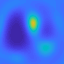}}
\put(190,360){\includegraphics[height=80pt]{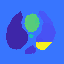}}

\put(10,270){\includegraphics[height=80pt]{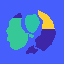}}
\put(100,270){\includegraphics[height=80pt]{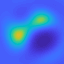}}
\put(190,270){\includegraphics[height=80pt]{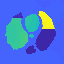}}

\put(10,180){\includegraphics[height=80pt]{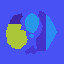}}
\put(100,180){\includegraphics[height=80pt]{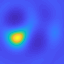}}
\put(190,180){\includegraphics[height=80pt]{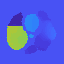}}


\put(10,90){\includegraphics[height=80pt]{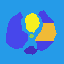}}
\put(100,90){\includegraphics[height=80pt]{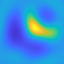}}
\put(190,90){\includegraphics[height=80pt]{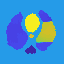}}

\put(10,0){\includegraphics[height=80pt]{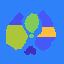}}
\put(100,0){\includegraphics[height=80pt]{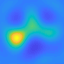}}
\put(190,0){\includegraphics[height=80pt]{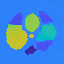}}

\put(38,455){\footnotesize \sc Truth}

\put(121,460){\footnotesize \sc Low Pass}
\put(115,450){\footnotesize \sc D-bar Image}

\put(203,460){\footnotesize \sc  Beltrami-Net}
\put(218,450){\footnotesize \sc  Image}

\put(0,30){\rotatebox{90}{{\sc \footnotesize Sim 5}}}
\put(0,120){\rotatebox{90}{{\sc \footnotesize Sim 4}}}
\put(0,210){\rotatebox{90}{{\sc \footnotesize Sim 3}}}
\put(0,300){\rotatebox{90}{{\sc \footnotesize Sim 2}}}
\put(0,390){\rotatebox{90}{{\sc \footnotesize Sim 1}}}

%
%
%
%
%
%
%
%
%

\end{picture}
\caption{Results for simulated test data with the network trained for the ACT4 data. Note that the training data only included single horizontal divisions in the lungs. Each row is plotted on its own scale.}  
\label{fig:ACT4demoSims}
\end{figure}

\clearpage

\subsection{Reconstructions from Experimental Data}\label{sec:resultsExp}
We next present reconstructions from the ACT4 and KIT4 experimental data.

\subsubsection{Experimental Reconstructions from ACT4}

\begin{figure}[h!]
\centering
\begin{picture}(300,370)

\put(10,255){\includegraphics[width=75pt]{HLSA_saline_DeepDbar_DICOM_smaller.jpg}}
\put(85,255){\includegraphics[width=75pt]{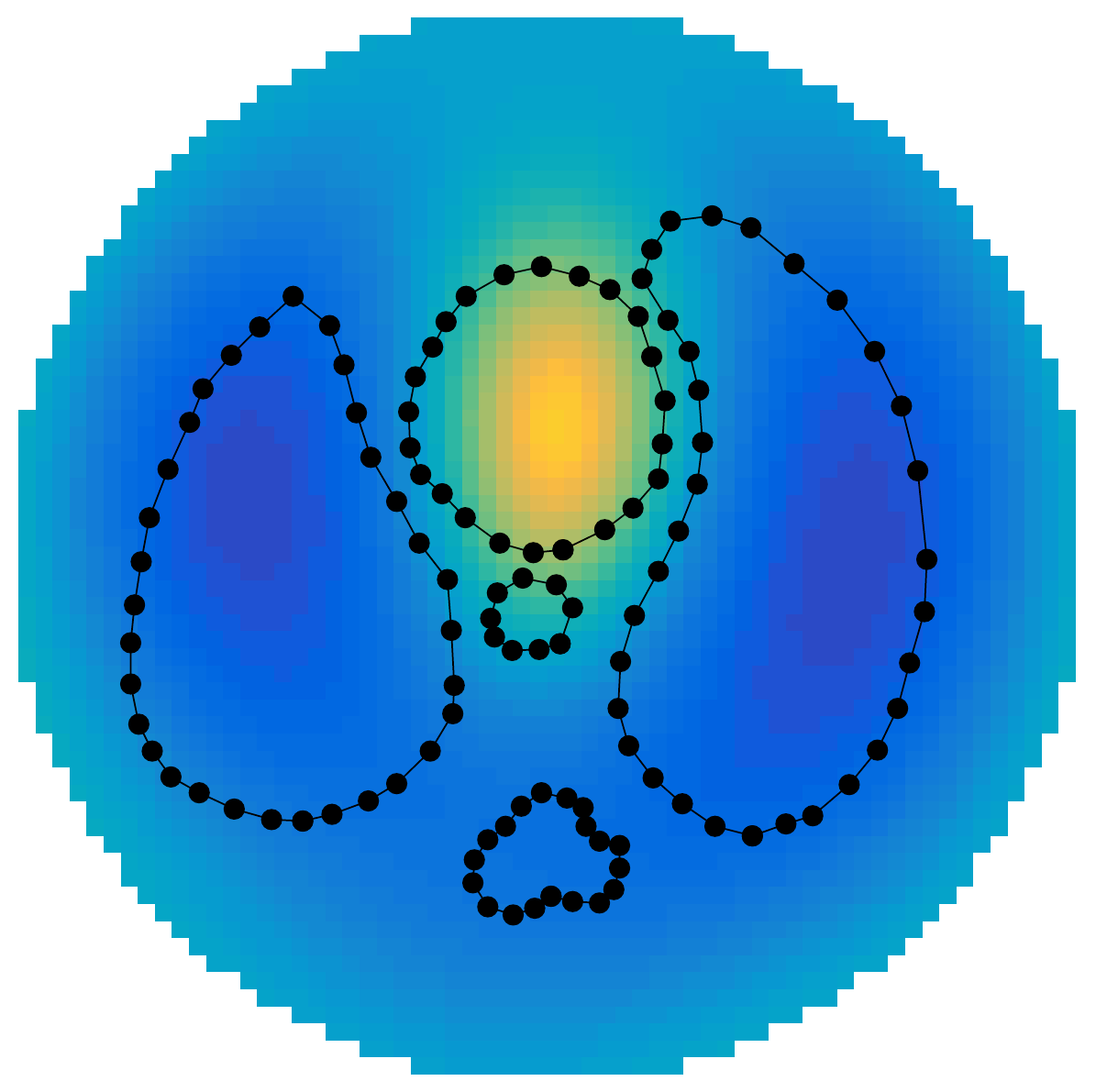}}
\put(160,255){\includegraphics[width=90pt]{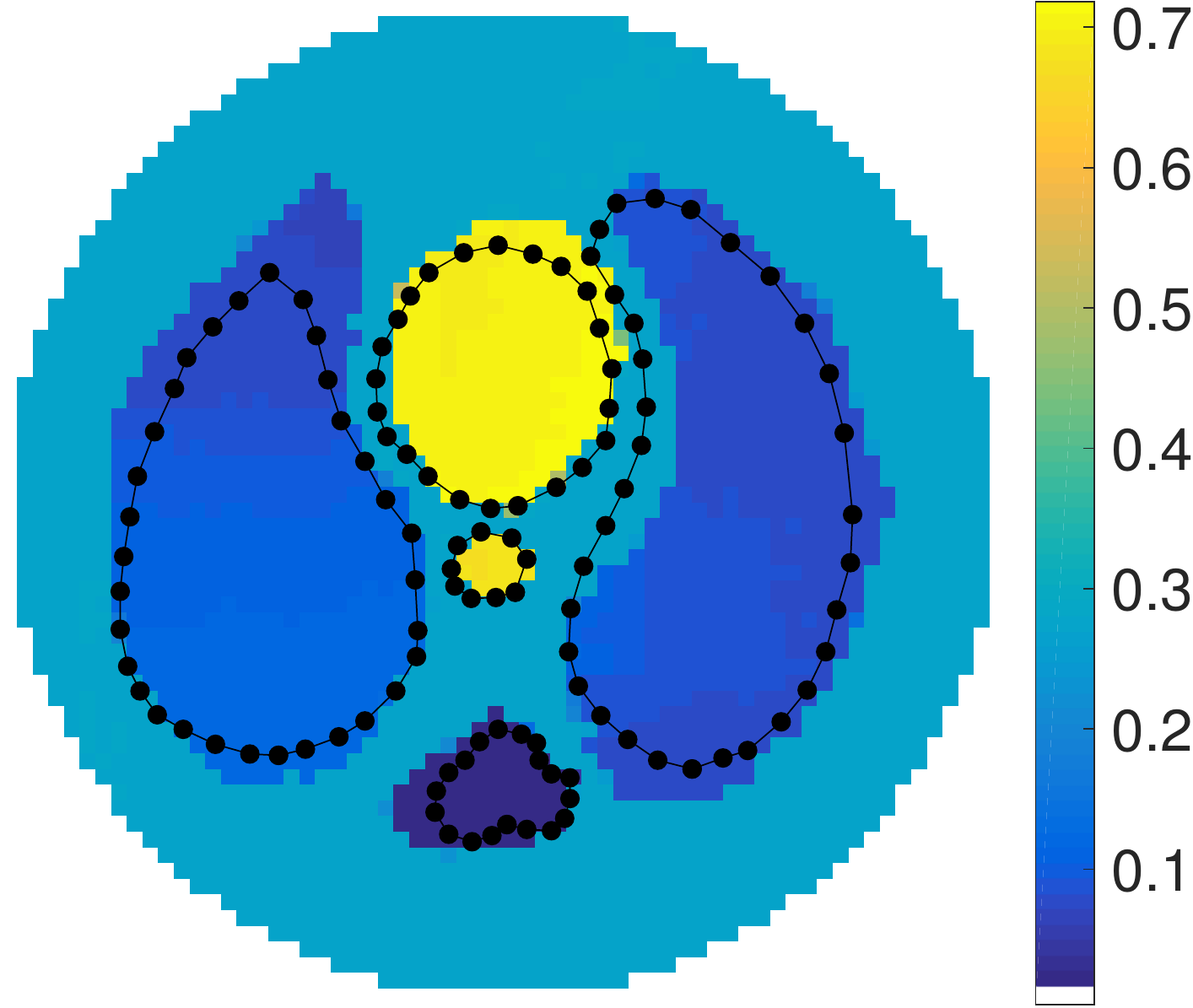}}

\put(10,170){\includegraphics[width=75pt]{HLSA_plueral_saline_DeepDbar_DICOM_smaller.jpg}}
\put(85,170){\includegraphics[width=75pt]{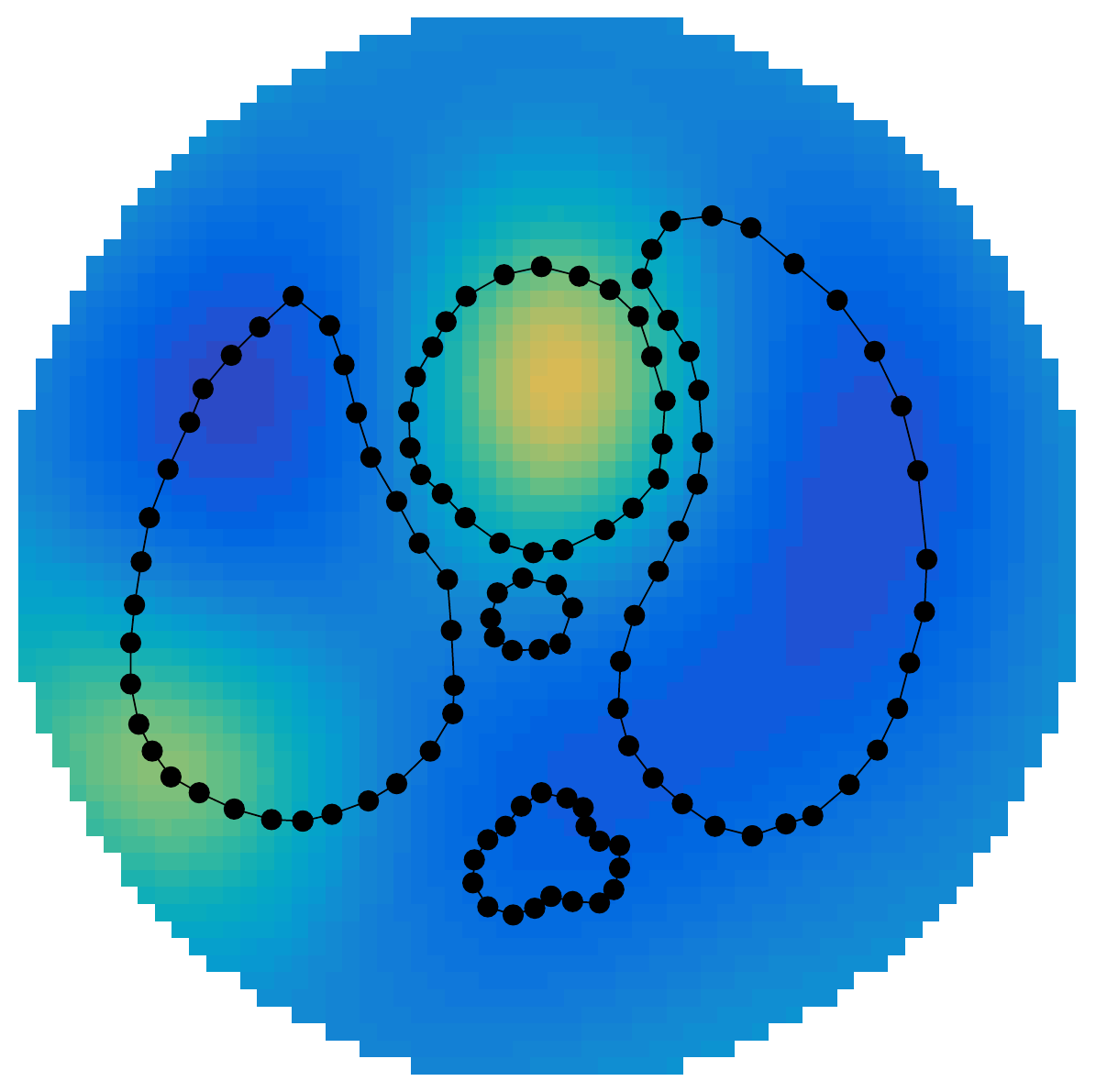}}
\put(160,170){\includegraphics[width=90pt]{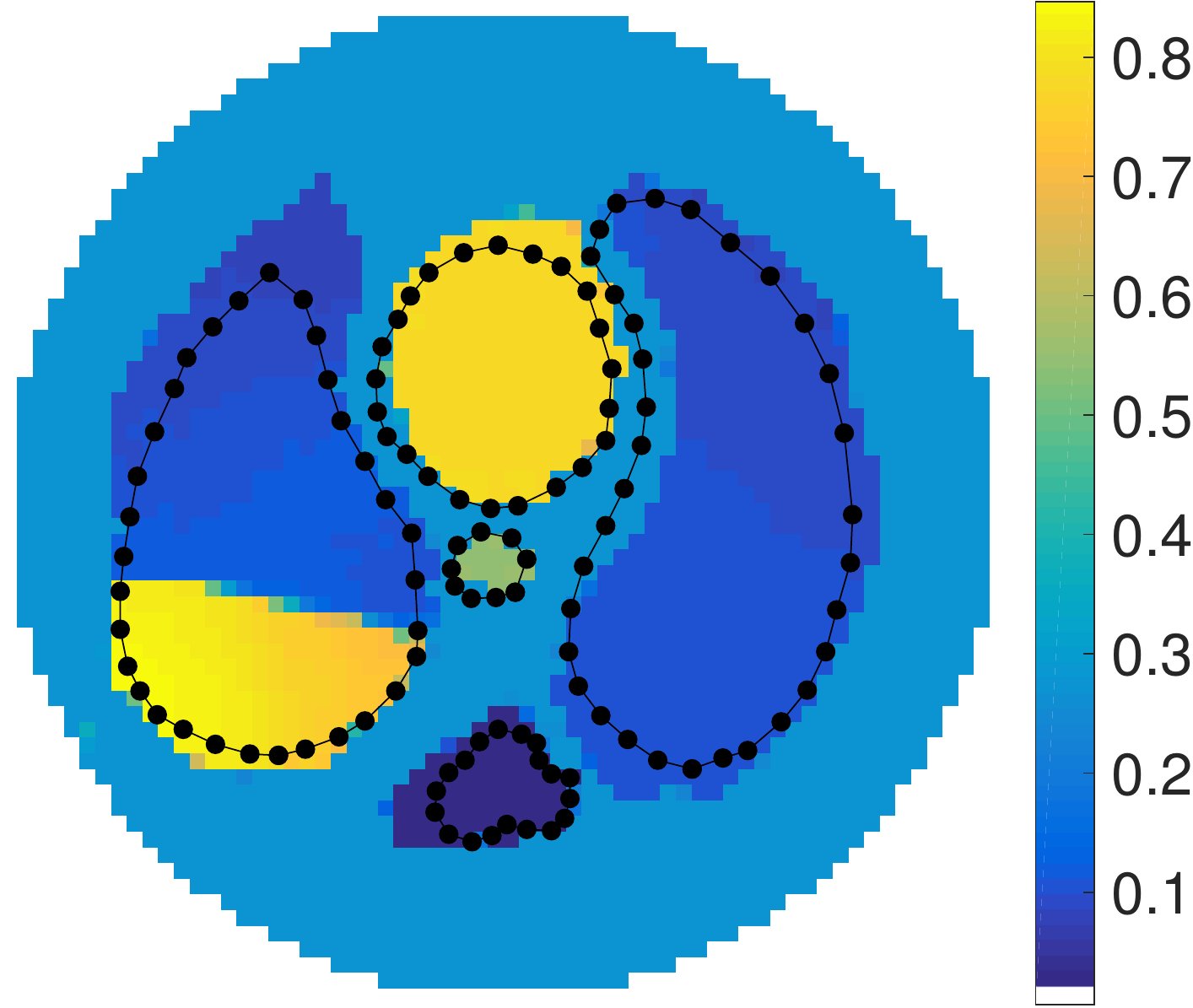}}

\put(10,85){\includegraphics[width=75pt]{HLSA_plastic_saline_DeepDbar_DICOM_smaller.jpg}}
\put(85,85){\includegraphics[width=75pt]{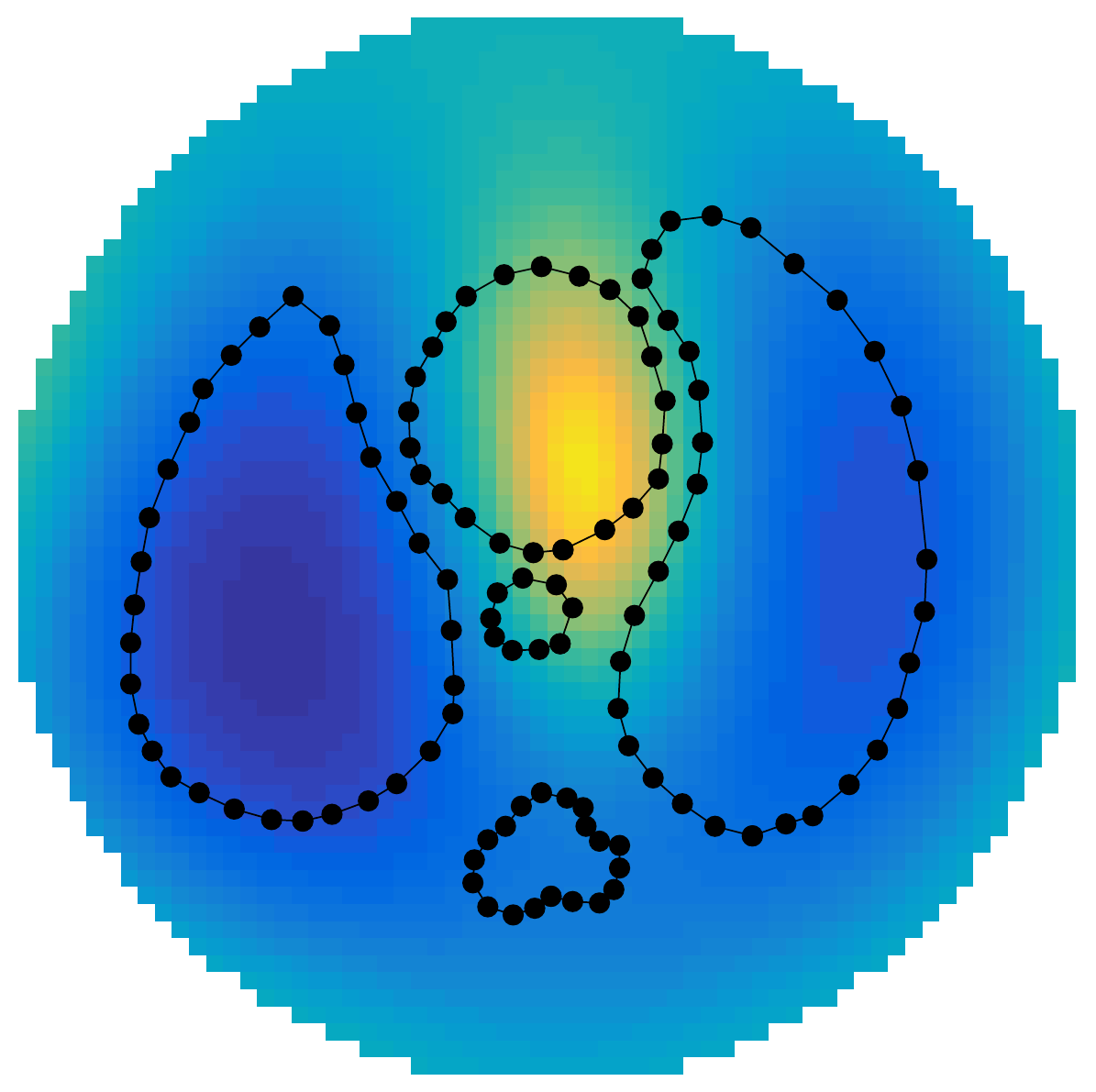}}
\put(160,85){\includegraphics[width=90pt]{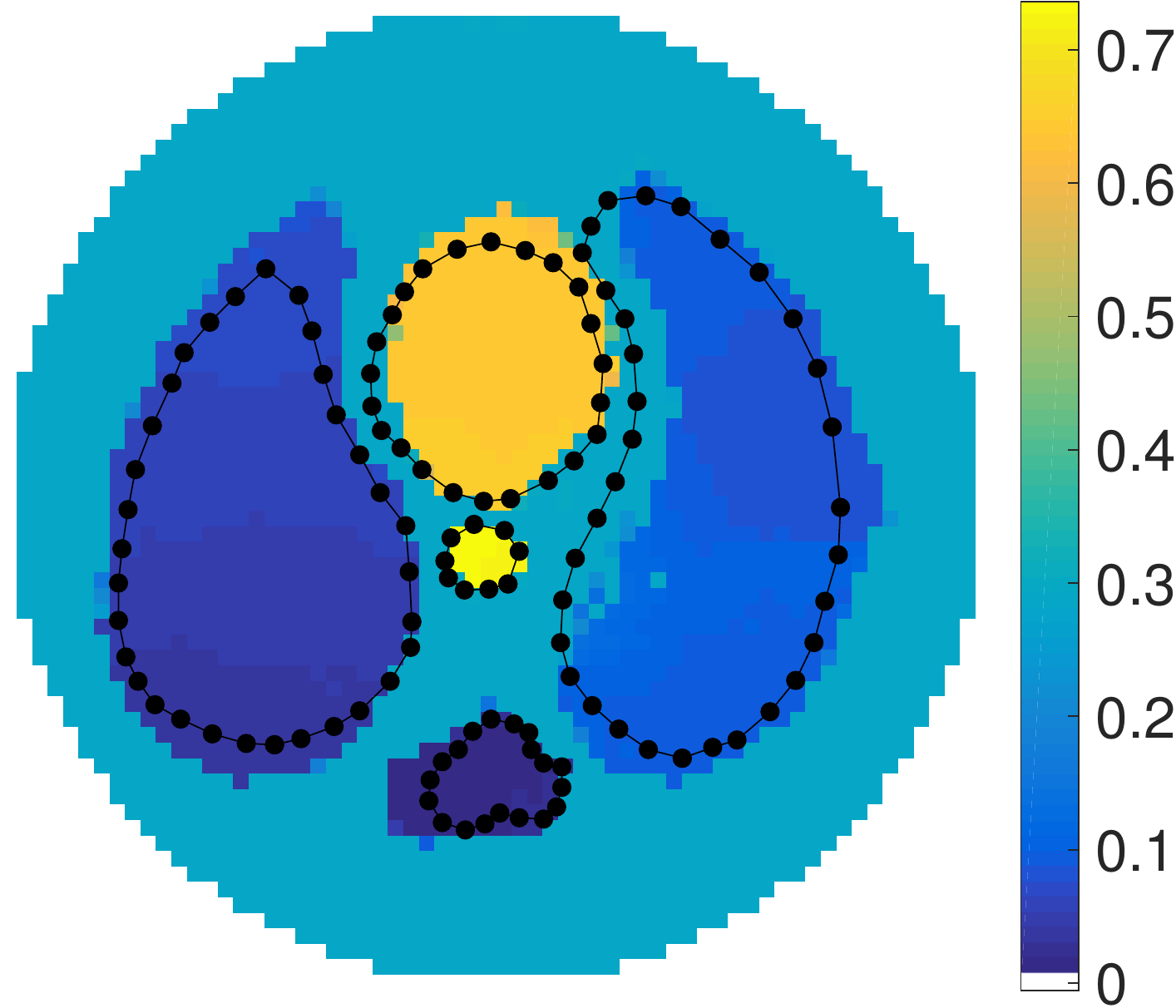}}

\put(10,0){\includegraphics[width=75pt]{ACT4_HLSA_wCopper_Correct_DICOM_smaller.jpg}}
\put(85,0){\includegraphics[width=75pt]{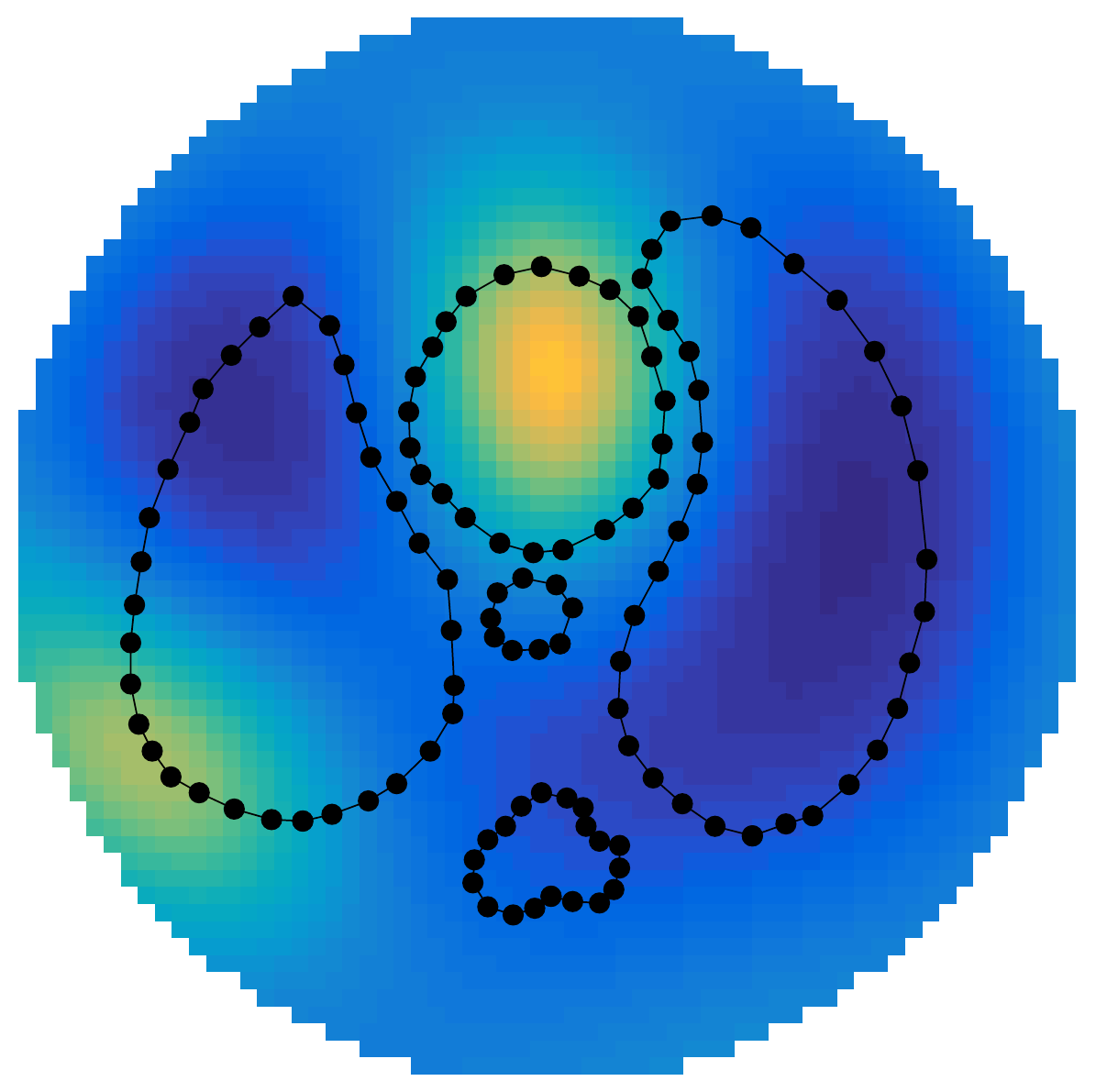}}
\put(160,0){\includegraphics[width=90pt]{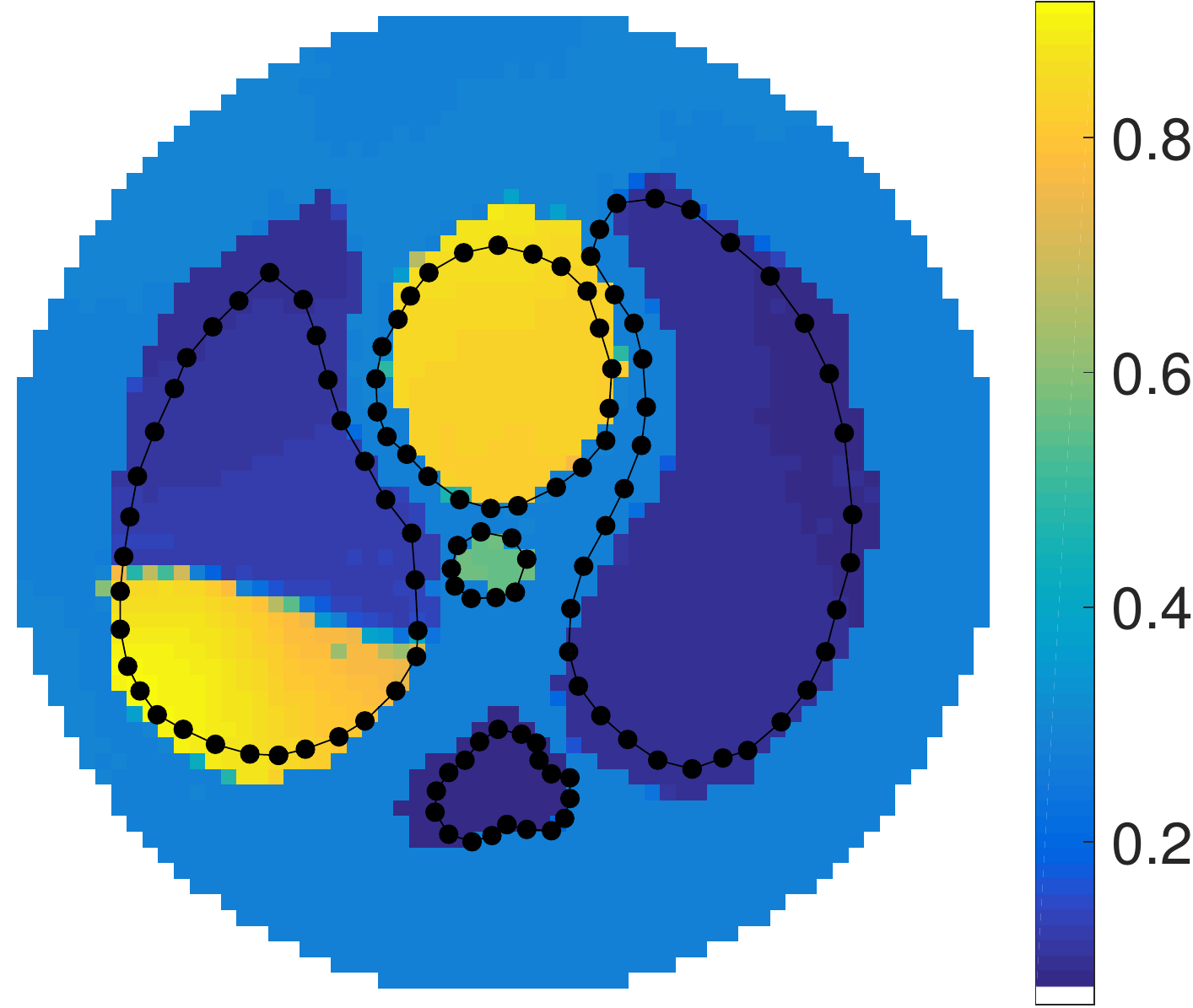}}

\put(25,350){\footnotesize \sc Experiment}

\put(108,355){\footnotesize \sc Low Pass}
\put(100,345){\footnotesize \sc D-bar Image}

\put(175,355){\footnotesize \sc  Beltrami-Net}
\put(187,345){\footnotesize \sc   Image}

\put(0,20){\rotatebox{90}{{\sc \footnotesize Injury 3}}}
\put(0,105){\rotatebox{90}{{\sc \footnotesize Injury 2}}}
\put(0,190){\rotatebox{90}{{\sc \footnotesize Injury 1}}}
\put(0,280){\rotatebox{90}{{\sc \footnotesize Healthy}}}

\end{picture}
\caption{\label{fig:ACT4_Results} Results for the experimental ACT4 data comparing the initial low-pass D-bar images to the post-processed Deep D-bar images.  Note that images are displayed here on the circular geometry of the tank, for presentation only.  The D-bar images on the full square $[-1,1]^2$ were used as inputs to the CNN to produce the Deep-Dbar images.  Each row is plotted on its own scale.}
\end{figure}

Figure~\ref{fig:ACT4_Results} depicts the results of the {\it Beltrami-Net} approach on four experiments with ACT4 data: {\sc Healthy} and {\sc Injuries 1-3} as shown in Figure~\ref{fig:RPI_phantoms}.  The black dots represent the approximate boundaries of the `healthy' organs, extracted from the photograph.  SSIMs, as well as relative $\ell_1$ and $\ell_2$ errors, were computed for the experimental reconstructions with the exception of {\sc Injury~3}, which has infinite conductors (copper tubes).  The comparisons, in Table~\ref{table:ACT4_quantResults}, used approximate `truth' images formed by assigning the measured conductivity values (Table~\ref{table:ACT4_setup}) in the respective regions.  Note that the coordinates for the bottom portion of the right (DICOM) lung were not specific to each injury, instead the entire region was assigned the same conductivity, even when the injury did not fill up the space as in Injury~2, plastic tubes and Injury~1 which is smaller than the original lung.

The obtained reconstructions for the ACT4 scenario are overall of high quality. Visually, we can identify the injuries in the lungs clearly from the Belrami-Net reconstructions as shown in Figure \ref{fig:ACT4_Results}. Both high conductive injuries are very clearly reconstructed and are even clearly visible in the D-Bar reconstructions. The lower conductive injury is harder to identify, in the D-bar reconstruction this results in a overall lower conductivity in the left lung. The Beltrami-Net then  manages to shift the lower conductivity to bottom of the lung, but can not establish a sharp boundary. We note here, that the network was only trained on horizontal injuries, nevertheless it manages to reproduce diagonal cuts for the high conductive injuries.

Quantitatively, the Beltrami-Net reconstructions show clear improvements over the low-pass D-bar reconstructions by all metrics in Table~\ref{table:ACT4_quantResults}. We remind here, that this is a case with strong a-priori knowledge and hence the results are expected to be of very high quality.  However, unlike the previous study, [\cite{Hamilton2018_DeepDbar}], the Beltrami-Net method did recover sharp diagonal divisions even when only training on horizontal cuts.


\begin{table}[h] 
\scriptsize
  \caption{Quantitative results for ACT4 experiments: Structural SIMilarity indices, as well as relative $\ell_1$ and $\ell_2$ images errors.}
    \begin{tabular}{l|c|c|c|c|c|c}
 & \multicolumn{3}{c}{\sc Low Pass D-Bar} & \multicolumn{3}{c}{\sc Beltrami-Net} \\
{\sc Experiment} &{\sc SSIM } &{\sc $\ell_1$-error} & {\sc $\ell_2$-error}  & {\sc SSIM } &{\sc $\ell_1$-error} & {\sc $\ell_2$-error} \\
    \hline
    \hline
    {\sc Healthy} 	& 0.5680 & 31.43\%	& 22.03\% & 0.7296 & 23.75\% & 13.75\% \\
    {\sc Agar} 		& 0.5176 & 35.87\%	& 24.62\% & 0.6963 & 27.79\% & 21.01\% \\
    {\sc Plastic} 	& 0.5085 & 34.91\%	& 24.44\% & 0.7053 & 22.26\% & 13.29\% \\
    \hline
    \hline
    \end{tabular}%
  \label{table:ACT4_quantResults}%
\end{table}%

\subsubsection{Experimental Reconstructions from KIT4}

\begin{figure}[h!]
\centering
\begin{picture}(300,360)

\put(10,255){\includegraphics[width=75pt]{KIT4_circTank_rotated_clipped.jpg}}
\put(85,255){\includegraphics[width=75pt]{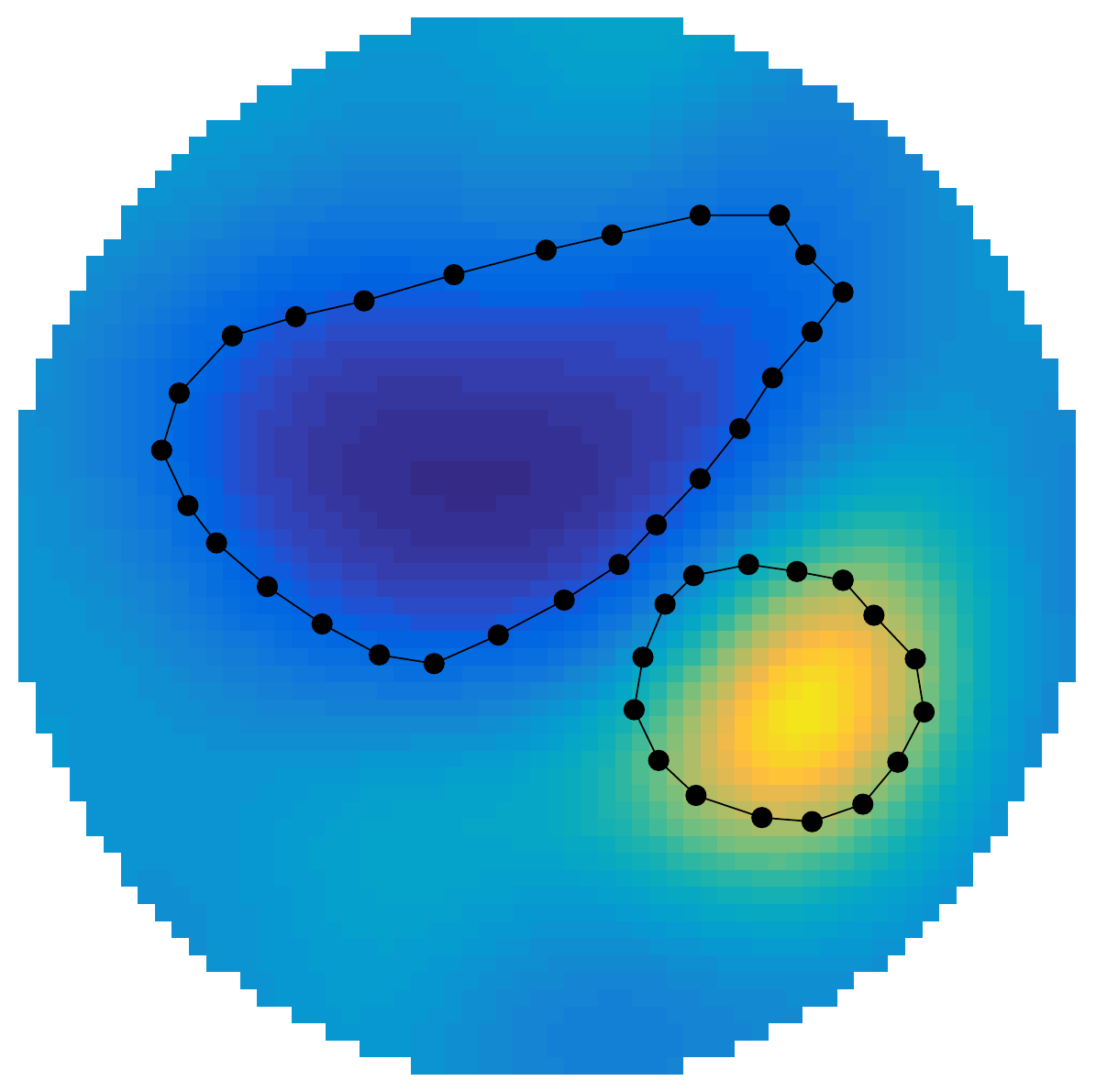}}
\put(160,255){\includegraphics[width=75pt]{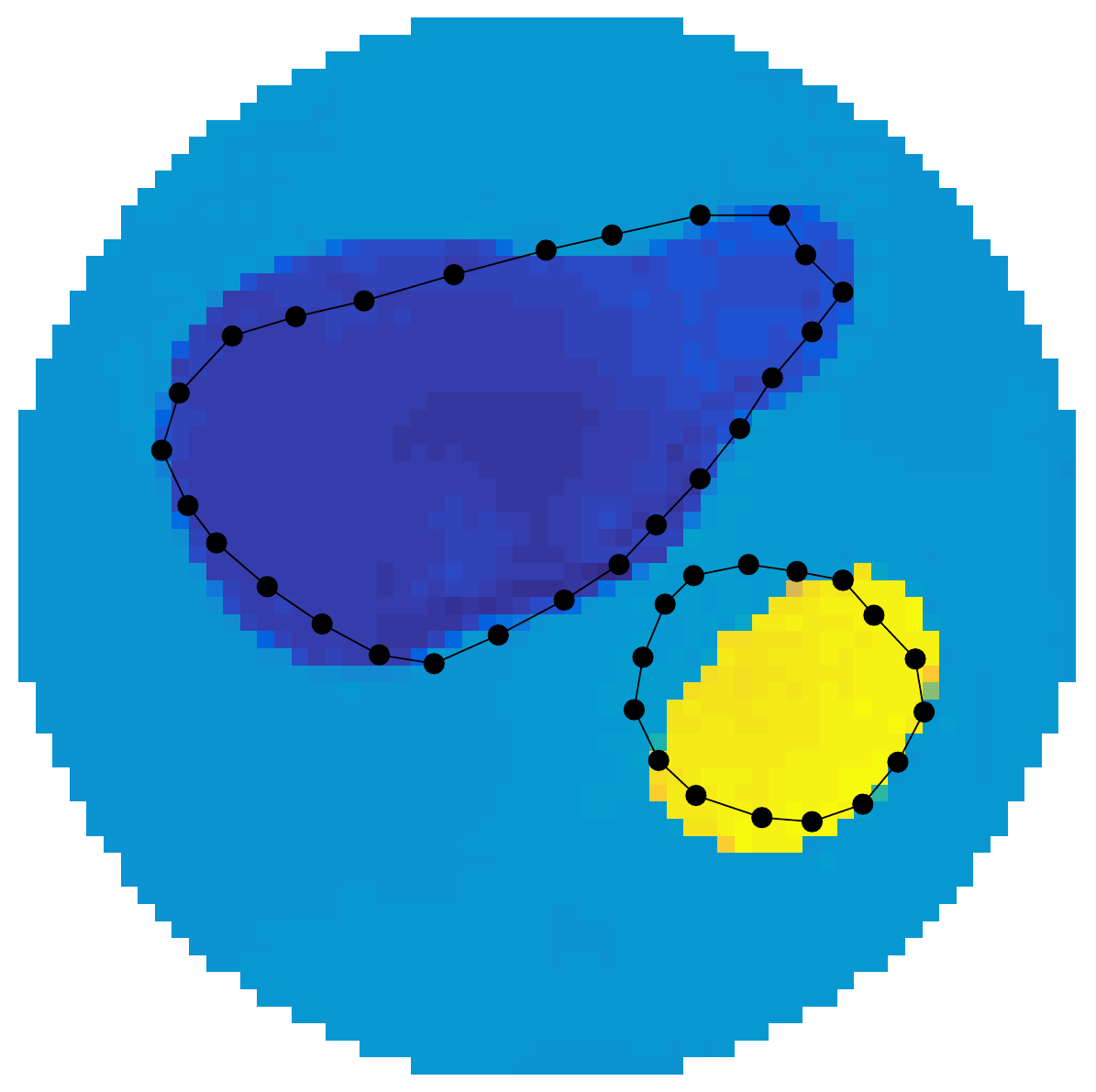}}
\put(237,255){\includegraphics[width=90pt]{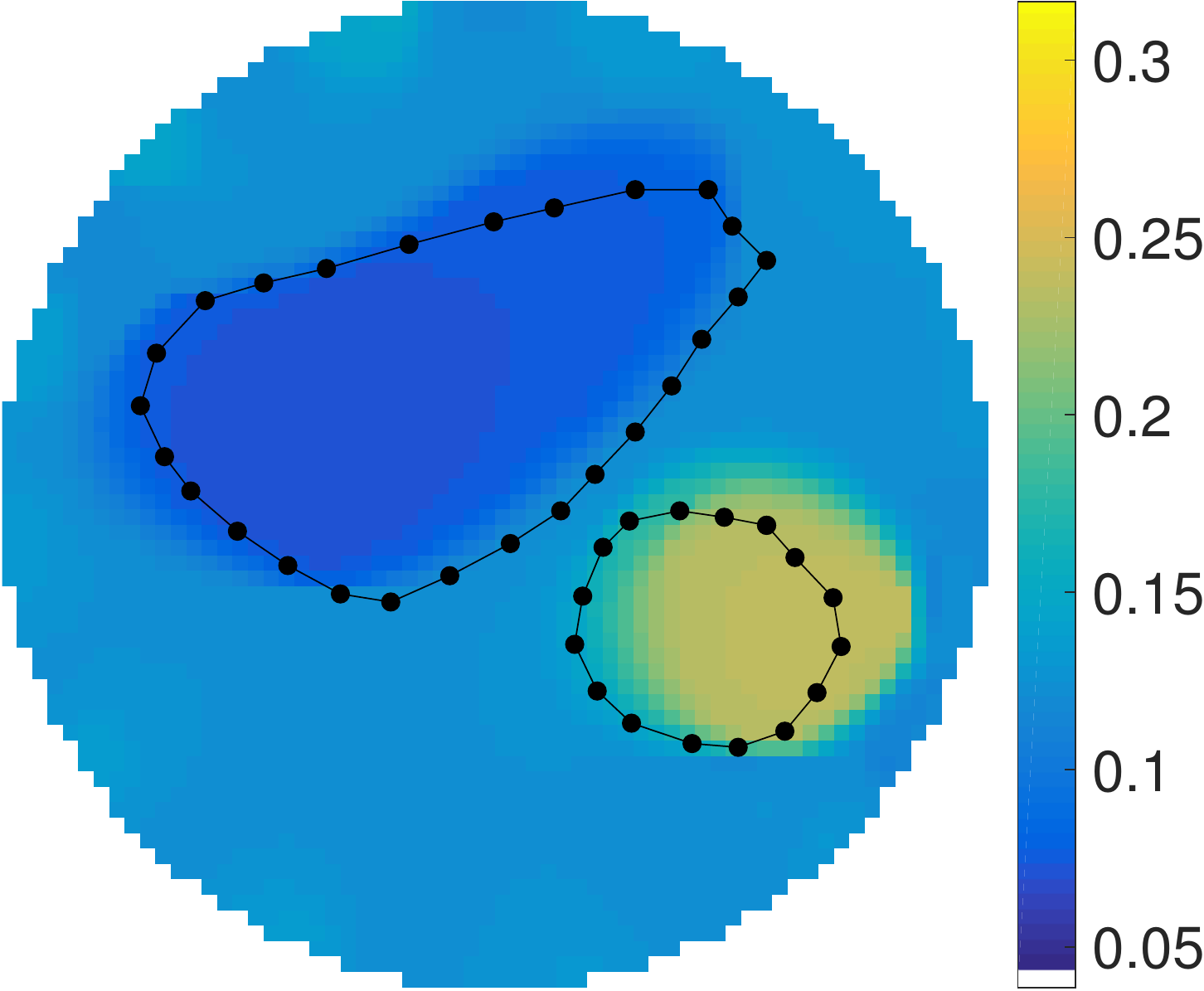}}

\put(10,170){\includegraphics[width=75pt]{KIT4_2018_chest_CORRECT.jpg}}
\put(85,170){\includegraphics[width=75pt]{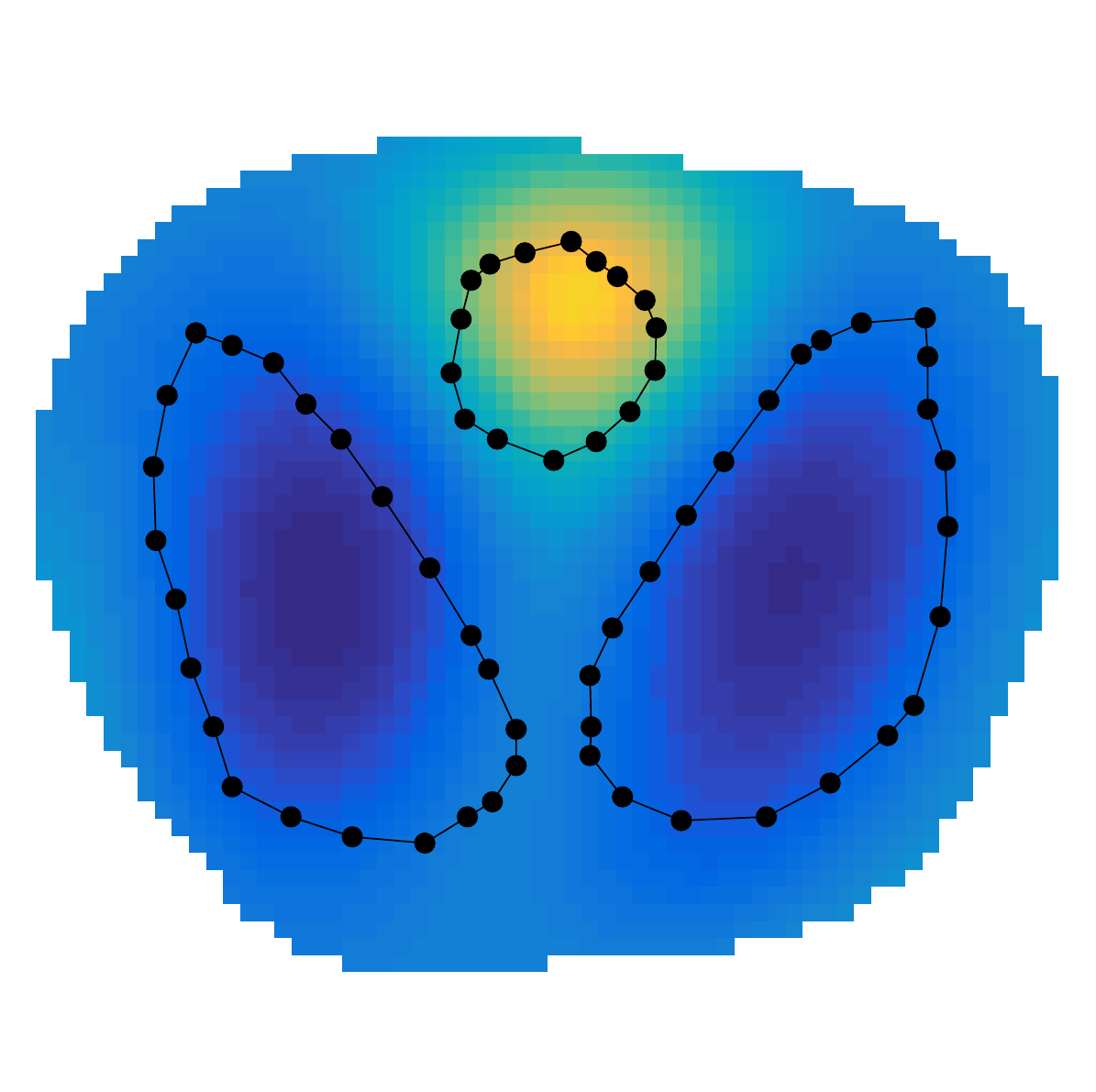}}
\put(160,170){\includegraphics[width=75pt]{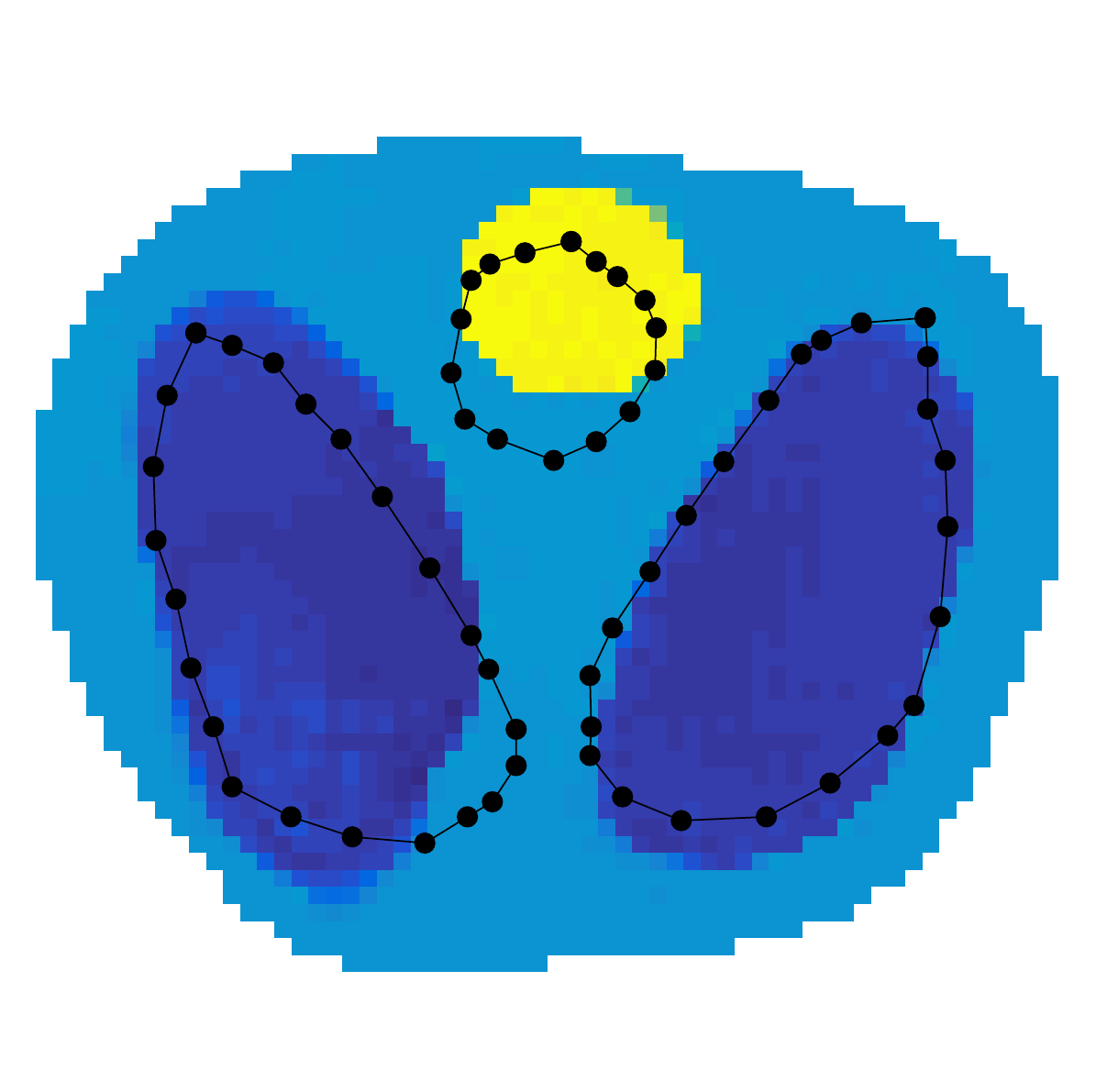}}
\put(237,170){\includegraphics[width=90pt]{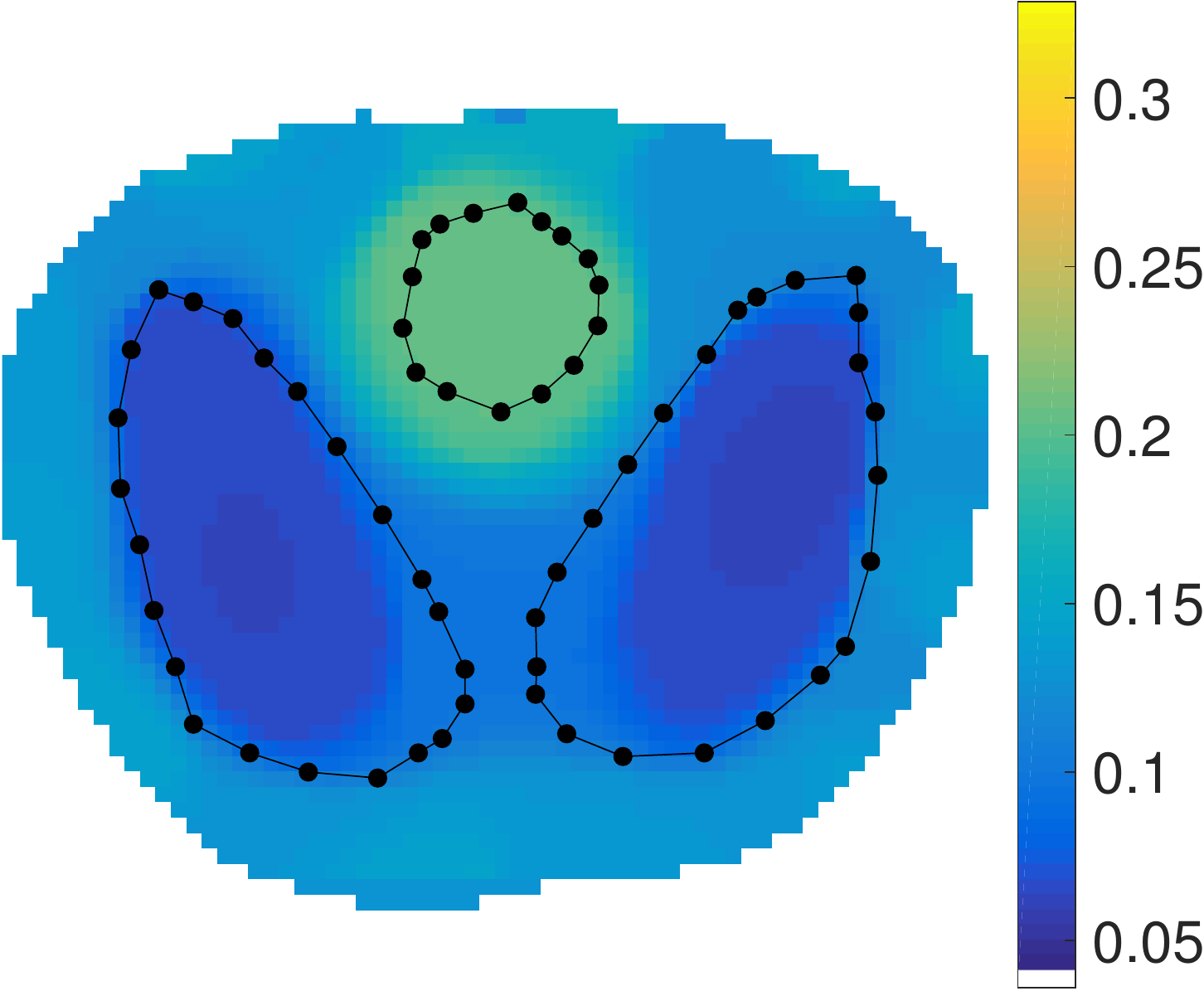}}

\put(10,85){\includegraphics[width=75pt]{KIT4_2018_chest_cut.jpg}}
\put(85,85){\includegraphics[width=75pt]{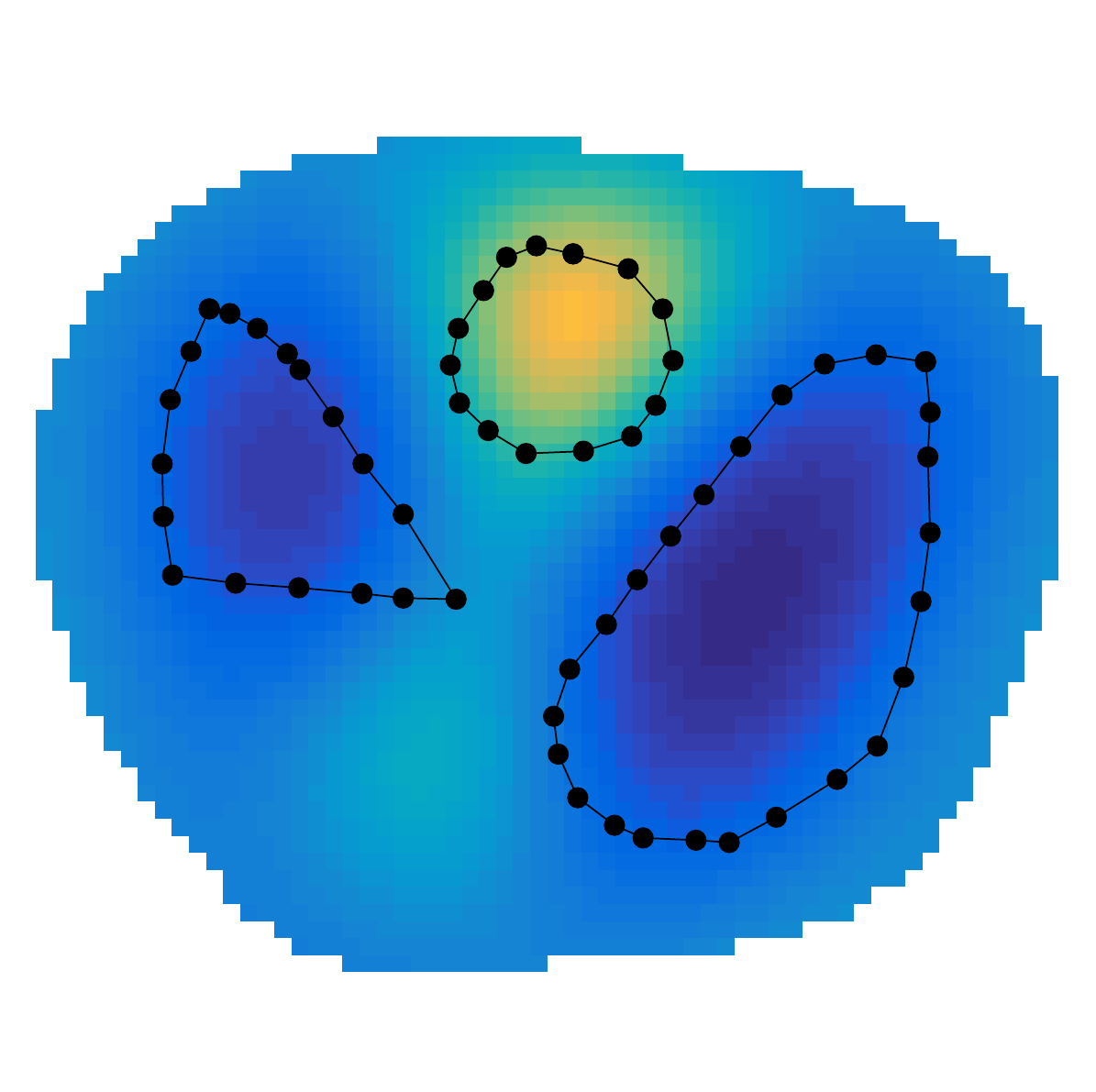}}
\put(160,85){\includegraphics[width=75pt]{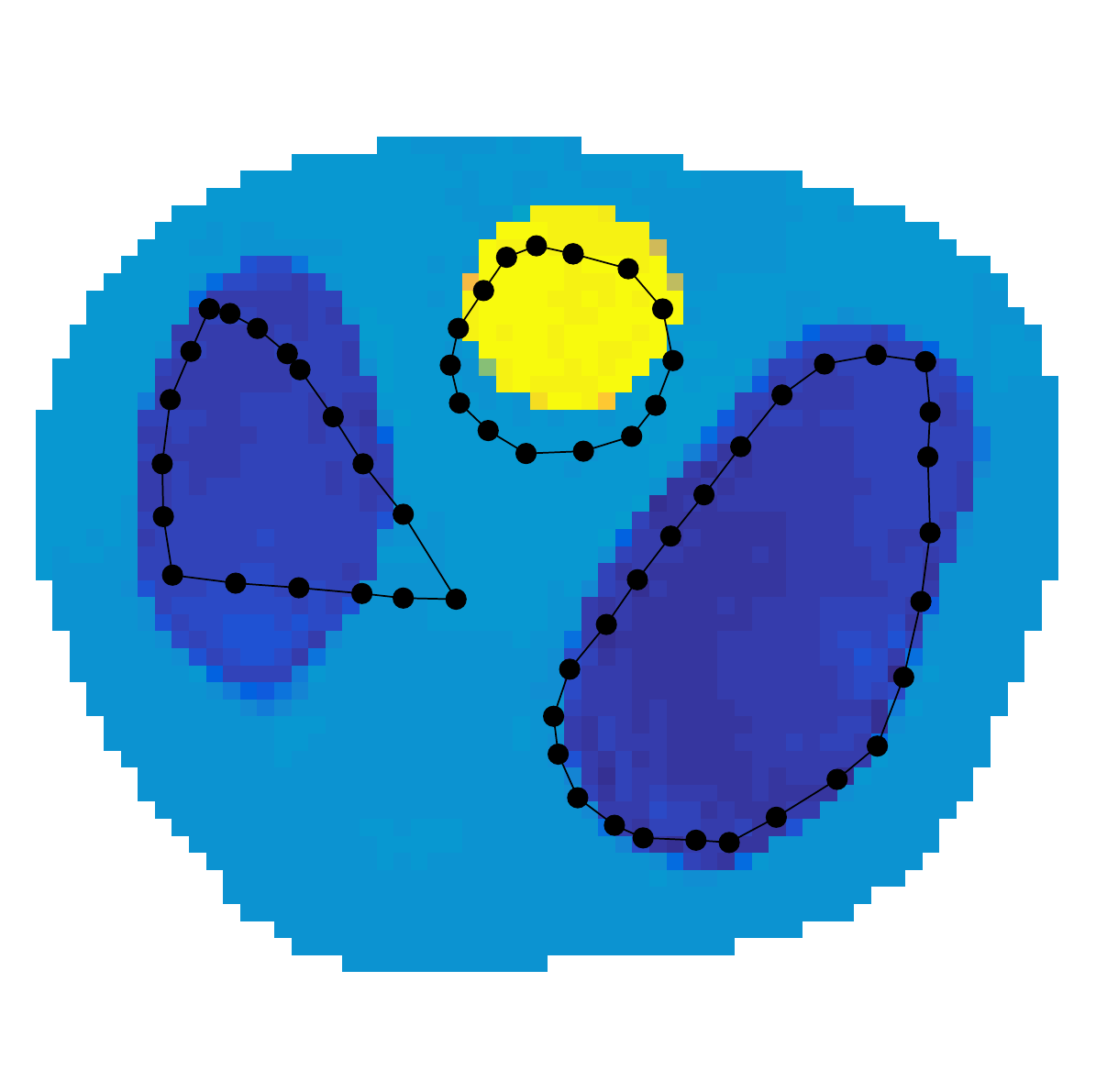}}
\put(237,85){\includegraphics[width=90pt]{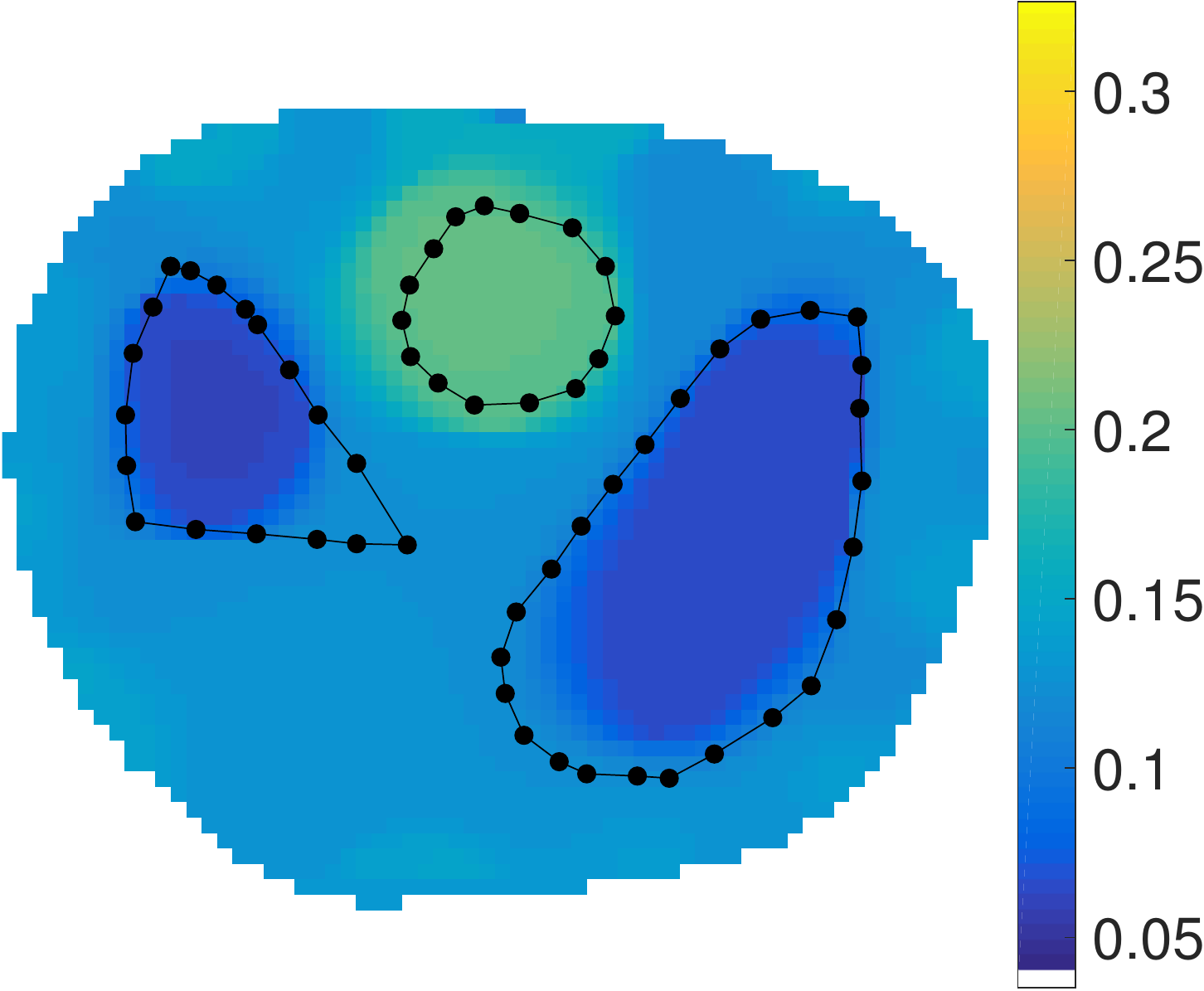}}

\put(10,0){\includegraphics[width=75pt]{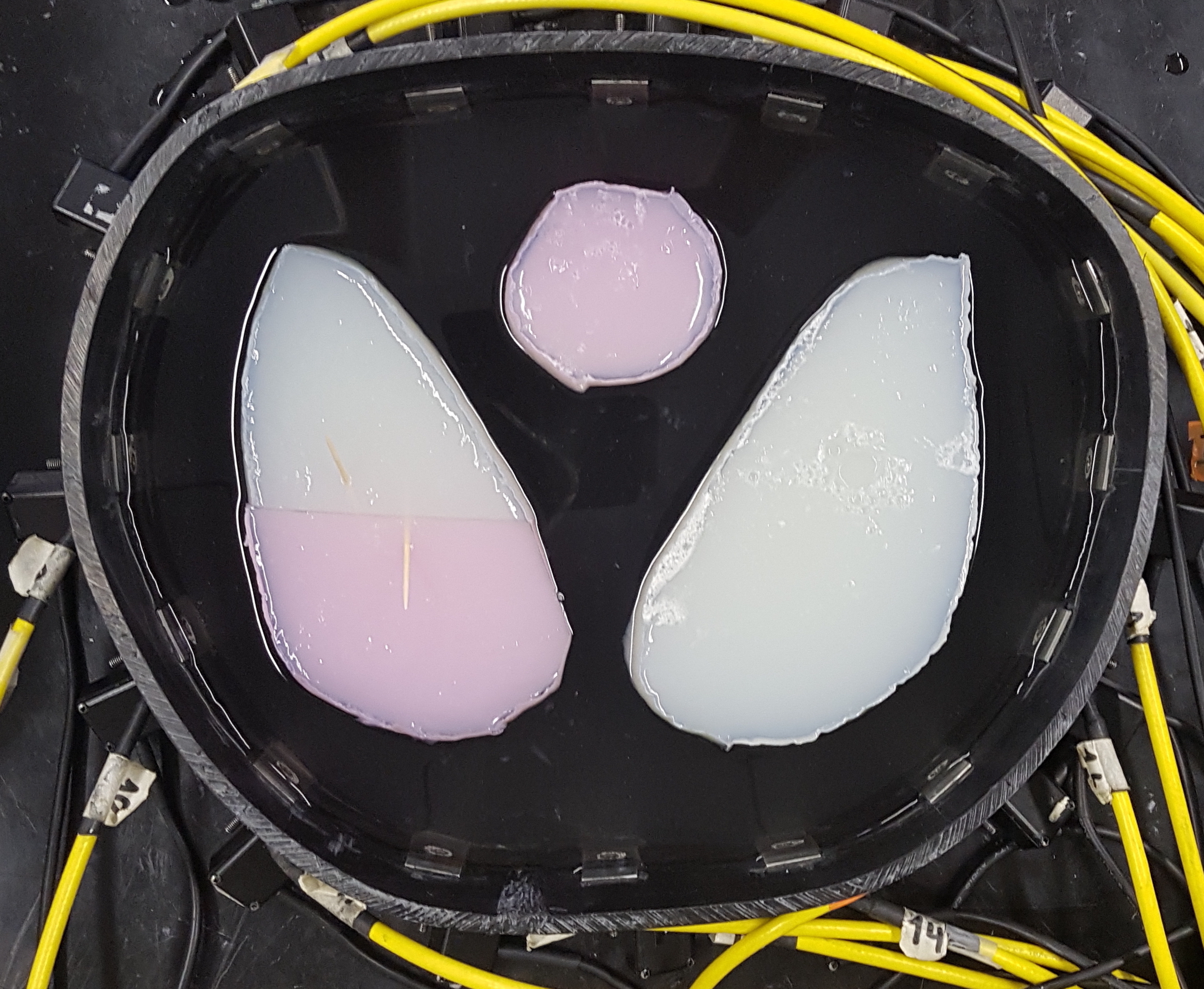}}
\put(85,0){\includegraphics[width=75pt]{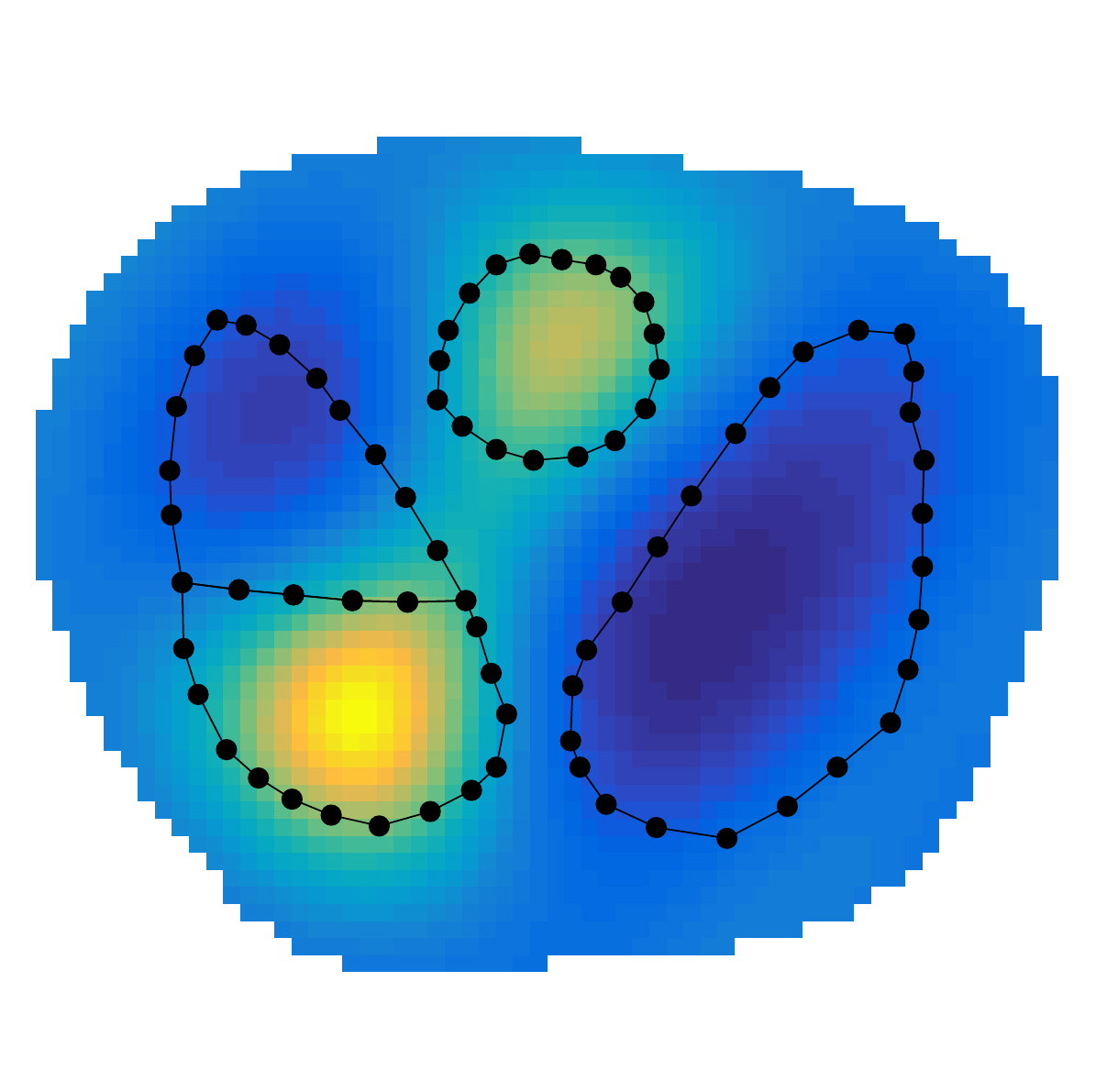}}
\put(160,0){\includegraphics[width=75pt]{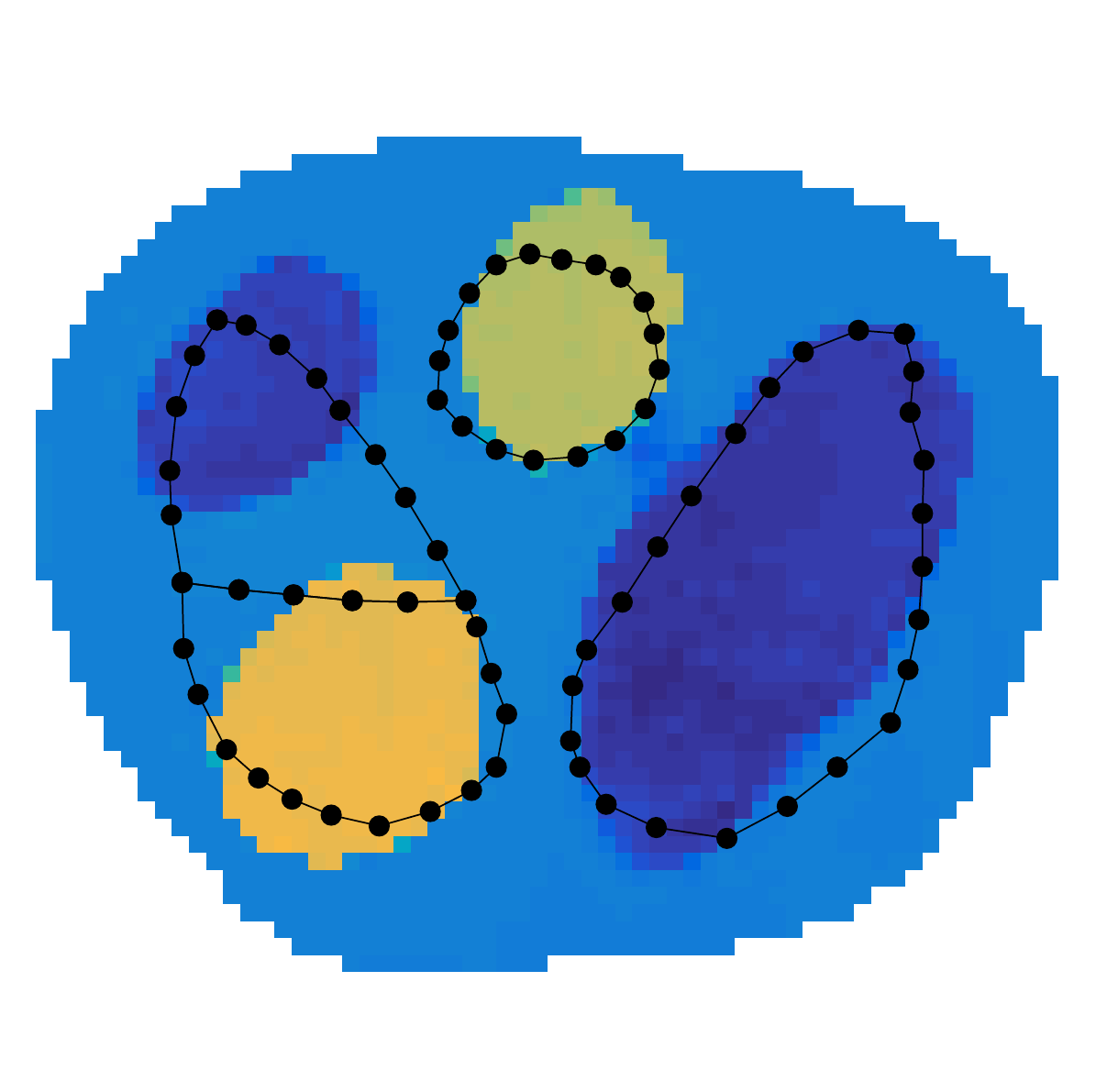}}
\put(237,0){\includegraphics[width=90pt]{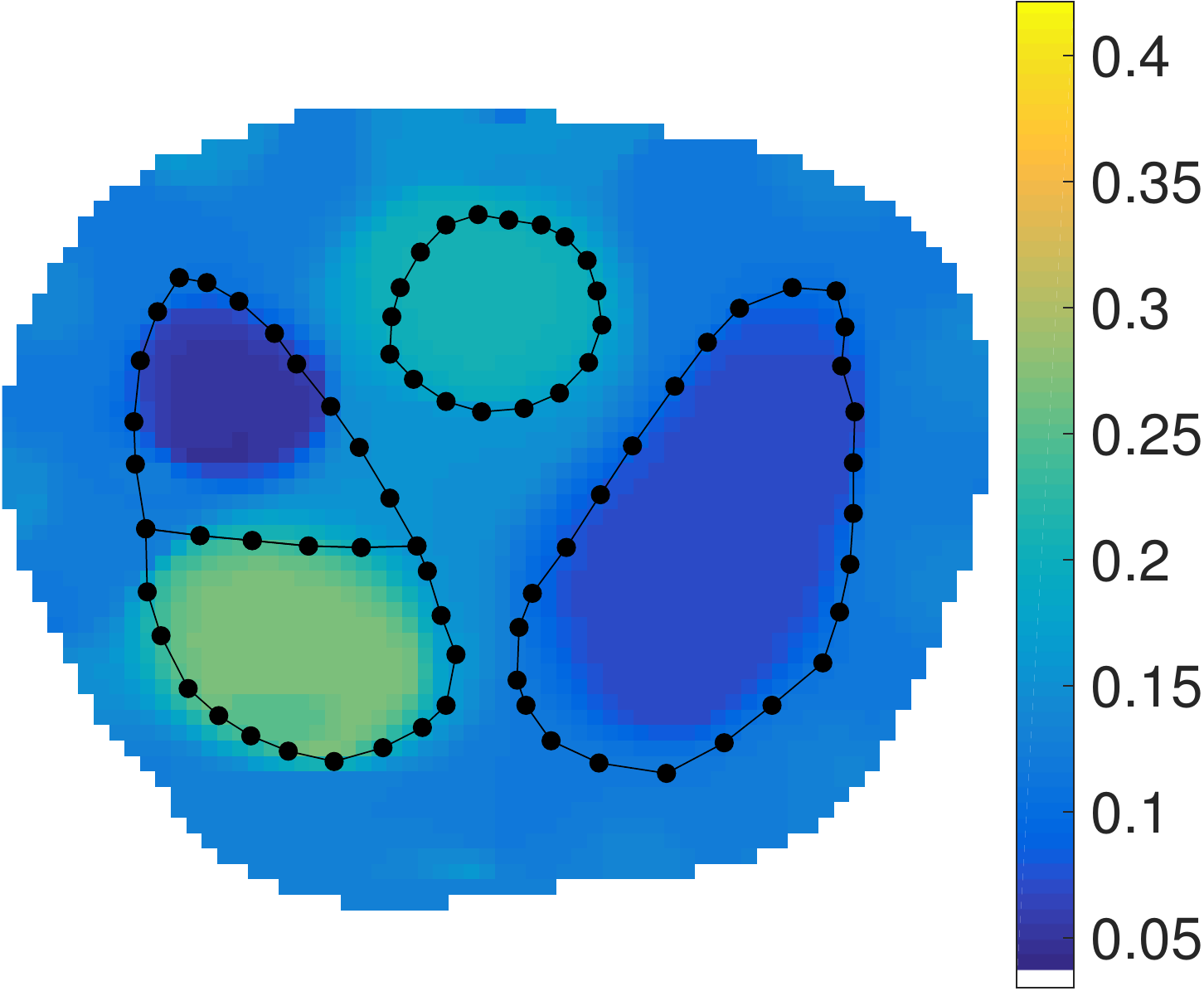}}

\put(25,350){\footnotesize \sc Experiment}

\put(108,355){\footnotesize \sc Low Pass}
\put(100,345){\footnotesize \sc D-bar Image}

\put(170,355){\footnotesize \sc  Beltrami-Net}
\put(185,345){\footnotesize \sc  Image}

\put(265,355){\footnotesize \sc  TV}
\put(260,345){\footnotesize \sc  Image}

\put(0,280){\rotatebox{90}{{\sc \footnotesize Circ}}}
\put(0,175){\rotatebox{90}{{\sc \footnotesize Chest-Healthy}}}
\put(0,95){\rotatebox{90}{{\sc \footnotesize Chest-Cut}}}
\put(0,10){\rotatebox{90}{{\sc \footnotesize Chest-Split}}}

\end{picture}
\caption{\label{fig:KIT4_Results2} KIT4 Results for the various test scenarios.  The initial D-bar image is compared to the Deep D-bar image.  The D-bar images, on the full square $[-1,1]^2$  are used as the `input' images for the CNN.  Images are displayed here clipped to their respective tank geometries for presentation only.  Each row is plotted on its own scale.}
\end{figure}

We next applied the {\it Beltrami-Net} method to the KIT4 datasets corresponding to Figure~\ref{fig:KIT4setups} and compared to total variation regularized least squares reconstructions (TV-LS) as outlined in Section~\ref{sec:Evaluation}. The reconstructed images are shown in Figure~\ref{fig:KIT4_Results2} and quantitative measurements (SSIM and relative $\ell_1$ and $\ell_2$ images errors) presented in Table~\ref{table:KIT4_quantResults}.

As one can see in Figure \ref{fig:KIT4_Results2}, all three methods produce images where the inclusions are clearly visible. The low-pass D-bar reconstructions are quite blurry as expected, but the post-processed images with Beltrami-Net are of very high contrast with sharp edges. In the TV-LS reconstructions, the boundary edges tend to be slightly blurred and there is a clear loss of contrast, which is a quite usual side-effect for TV regularized reconstructions. Neither of the methods is able to identify the split chest in the fourth phantom, and instead separate the lung into the two areas of opposing conductivity with saline between them. We note here that the Beltrami-Net was trained with minimal prior knowledge of only elliptic inclusions. Nevertheless, the Beltrami-Net reconstructions show shapes that differ from this simple prior. Hence we hypothesize that the network mainly learns a segmentation and correction of the existing features in the D-bar reconstructions.

The quantitative measures, SSIM, as well as relative $\ell_1$ and $\ell_2$ image errors, were computed for each case by comparing to approximate `truth' images constructed using the measured conductivity values and photographs of the experiments, see Table~\ref{table:KIT4_quantResults}. The quantitative improvements of Beltrami-Net are rather minor in this case. This is as expected due to low prior information. SSIM of D-Bar and Beltrami-Net are quite comparable, but generally high already. Most notably, even though the $\ell^2$-error is quite constant as well, there is a clear improvement in $\ell^1$-error, most likely due to sharper boundary edges.  The TV-LS method provides comparable metrics and reconstructions, outperforming both the low-pass D-bar and Beltrami-Net methods for the SSIM of the {\it Chest Healthy} and {\it Chest Cut} phantoms, but underperforming for the {\it Chest Split} experiment. Most notably, the Beltrami-Net reconstruction are consistently better in $\ell^1$-error for all provided measures.

\begin{table}[h] 
\scriptsize
  \caption{Quantitative results for KIT4 experiments.}
    \begin{tabular}{l|c|c|c|c|c|c|c|c|c}
 & \multicolumn{3}{c}{\sc Low Pass D-Bar} & \multicolumn{3}{c}{\sc Beltrami-Net}& \multicolumn{3}{c}{\sc TV}  \\
{\sc Experiment} &{\sc SSIM } &{\sc $\ell_1$-error} & {\sc $\ell_2$-error}  & {\sc SSIM } &{\sc $\ell_1$-error} & {\sc $\ell_2$-error} & {\sc SSIM } &{\sc $\ell_1$-error} & {\sc $\ell_2$-error}\\
    \hline
    \hline
    {\sc Circ Agar} 	& 0.8831 & 23.08\% & 14.39\% & 0.8921 &	19.53\% & 13.11\% &0.8843	&22.09\%	&16.14\%\\
    {\sc Chest Healthy} & 0.8507 & 26.29\% & 15.73\% & 0.8370 & 21.03\% & 17.33\% &0.8709 &	24.30\%&	17.03\%\\
    {\sc Chest Cut} 	& 0.8684 & 22.56\% & 15.55\% & 0.8516 & 18.67\% & 15.26\% &0.8939&	20.80\% &16.11\%\\
    {\sc Chest Split}	& 0.8244 & 28.79\% & 14.76\% & 0.8267 &	21.78\%	& 16.90\% &0.7877	 &36.28\%	&36.25\% \\
    \hline
    \hline
    \end{tabular}%
  \label{table:KIT4_quantResults}%
\end{table}%
\section{Conclusions}\label{sec:conc}
In this work we considered two conceptually different settings: i) A constrained case of thoracic imaging with the ACT4 measurements, where high a-priori knowledge is available, and ii) A very general setting with the KIT4 experiments on varying tank boundary and inclusion shapes with minimal prior knowledge in the training data. Consequently, the obtained results are slightly different in their nature. Whereas the ACT4 reconstructions are of very high quality and close to the target/image prior, the KIT4 reconstructions are more general and it is harder to obtain the exact shapes of the targets. Compared the the reference method of total variation constrained least square reconstructions, the reconstruction quality of Beltrami-Net is quite similar with a slight advantage in contrast and hence $\ell^1$-error measures.

We believe that this comparison provides good insight of what is possible in EIT in combination with deep learning based post-processing, in particular for D-bar based methods. We remind here, that EIT is a highly ill-posed inverse problem and hence it is not surprising that strong prior knowledge is needed to obtain high-quality images. Thus, we believe that the presented approach will be most useful in constrained imaging settings, where boundary shapes might vary, such as thoracic imaging for the identification of lung volumes or injuries. Additionally, process monitoring and non-destructive testing, where knowledge of possible composition and defects is known, may be areas of interest for this approach.

\section*{Acknowledgements}
We acknowledge the support of NVIDIA Corporation for the donation of the Titan Xp GPU used for this research.  A.~H\"anninen and V.~Kolehmainen acknowledge the Academy of Finland (Project 312343, Finnish Centre of Excellence in Inverse Modelling and Imaging) and the Jane and Aatos Erkko Foundation.

\small
\bibliography{bibliographyRefs_Ville}

\begin{thebibliography}{33}
\providecommand{\natexlab}[1]{#1}
\providecommand{\url}[1]{\texttt{#1}}
\expandafter\ifx\csname urlstyle\endcsname\relax
  \providecommand{\doi}[1]{doi: #1}\else
  \providecommand{\doi}{doi: \begingroup \urlstyle{rm}\Url}\fi

\bibitem[Adler and {\"O}ktem(2017)]{Adler2017}
Jonas Adler and Ozan {\"O}ktem.
\newblock Solving ill-posed inverse problems using iterative deep neural
  networks.
\newblock \emph{Inverse Problems}, 33\penalty0 (12):\penalty0 124007, 2017.

\bibitem[Astala and P{\"a}iv{\"a}rinta(2006{\natexlab{a}})]{Astala2006}
K.~Astala and L.~P{\"a}iv{\"a}rinta.
\newblock A boundary integral equation for {C}alder\'on's inverse conductivity
  problem.
\newblock In \emph{Proc. 7th Internat. Conference on Harmonic Analysis,
  Collectanea Mathematica}, 2006{\natexlab{a}}.

\bibitem[Astala and P{\"a}iv{\"a}rinta(2006{\natexlab{b}})]{Astala2006a}
K.~Astala and L.~P{\"a}iv{\"a}rinta.
\newblock {C}alder\'on's inverse conductivity problem in the plane.
\newblock \emph{Annals of Mathematics}, 163\penalty0 (1):\penalty0 265--299,
  2006{\natexlab{b}}.
\newblock ISSN 0003-486X.
\newblock \doi{10.4007/annals.2006.163.265}.
\newblock URL \url{http://dx.doi.org/10.4007/annals.2006.163.265}.

\bibitem[Astala et~al.(2010)Astala, Mueller, P{\"a}iv{\"a}rinta, and
  Siltanen]{Astala2010}
K.~Astala, J.L. Mueller, L.~P{\"a}iv{\"a}rinta, and S.~Siltanen.
\newblock Numerical computation of complex geometrical optics solutions to the
  conductivity equation.
\newblock \emph{Applied and Computational Harmonic Analysis}, 29\penalty0
  (1):\penalty0 391--403, 2010.

\bibitem[Cheney et~al.(1990)Cheney, Isaacson, and Isaacson]{Cheney1990}
M.~Cheney, D.~Isaacson, and E.L. Isaacson.
\newblock Exact solutions to a linearized inverse boundary value problem.
\newblock \emph{Inverse Problems}, 6:\penalty0 923--934, 1990.

\bibitem[Cheney et~al.(1999)Cheney, Isaacson, and Newell]{Cheney1999}
M.~Cheney, D.~Isaacson, and J.~C. Newell.
\newblock Electrical impedance tomography.
\newblock \emph{SIAM Review}, 41\penalty0 (1):\penalty0 85--101, 1999.

\bibitem[Dodd and Mueller(2014)]{Dodd2014}
M.~Dodd and J.~L. Mueller.
\newblock A real-time {D}-bar algorithm for 2-{D} electrical impedance
  tomography data.
\newblock \emph{Inverse Problems and Imaging}, 8\penalty0 (4):\penalty0
  1013--1031, 2014.
\newblock ISSN 1930-8337.
\newblock \doi{10.3934/ipi.2014.8.1013}.
\newblock URL
  \url{http://aimsciences.org/journals/displayArticlesnew.jsp?paperID=10528}.

\bibitem[González et~al.(2017)González, Kolehmainen, and
  Seppänen]{Gonzalez2017}
Gerardo González, Ville Kolehmainen, and Aku Seppänen.
\newblock Isotropic and anisotropic total variation regularization in
  electrical impedance tomography.
\newblock \emph{Computers \& Mathematics with Applications}, 74\penalty0
  (3):\penalty0 564 -- 576, 2017.
\newblock ISSN 0898-1221.
\newblock \doi{https://doi.org/10.1016/j.camwa.2017.05.004}.
\newblock URL
  \url{http://www.sciencedirect.com/science/article/pii/S0898122117302833}.

\bibitem[Hamilton and Hauptmann(2018)]{Hamilton2018_DeepDbar}
S.~J. Hamilton and A.~Hauptmann.
\newblock Deep d-bar: Real time electrical impedance tomography imaging with
  deep neural networks.
\newblock \emph{IEEE Transactions on Medical Imaging}, 37\penalty0
  (10):\penalty0 2367--2377, 2018.

\bibitem[Hamilton et~al.(2018)Hamilton, Mueller, and
  Santos]{Hamilton2018_Robust}
S~J Hamilton, J~L Mueller, and T~R Santos.
\newblock Robust computation in 2d absolute eit (a-eit) using d-bar methods
  with the ‘exp’ approximation.
\newblock \emph{Physiological Measurement}, 39\penalty0 (6):\penalty0 064005,
  2018.
\newblock URL \url{http://stacks.iop.org/0967-3334/39/i=6/a=064005}.

\bibitem[Hammernik et~al.(2018)Hammernik, Klatzer, Kobler, Recht, Sodickson,
  Pock, and Knoll]{Hammernik2018}
Kerstin Hammernik, Teresa Klatzer, Erich Kobler, Michael~P Recht, Daniel~K
  Sodickson, Thomas Pock, and Florian Knoll.
\newblock Learning a variational network for reconstruction of accelerated mri
  data.
\newblock \emph{Magnetic resonance in medicine}, 79\penalty0 (6):\penalty0
  3055--3071, 2018.

\bibitem[Hauptmann et~al.(2018{\natexlab{a}})Hauptmann, Arridge, Lucka,
  Muthurangu, and Steeden]{Hauptmann2018}
A~Hauptmann, S~Arridge, F~Lucka, V~Muthurangu, and JA~Steeden.
\newblock Real-time cardiovascular mr with spatio-temporal artifact suppression
  using deep learning-proof of concept in congenital heart disease.
\newblock \emph{Magnetic resonance in medicine}, 2018{\natexlab{a}}.

\bibitem[Hauptmann(2017)]{Hauptmann2017}
Andreas Hauptmann.
\newblock Approximation of full-boundary data from partial-boundary electrode
  measurements.
\newblock \emph{Inverse Problems}, 33\penalty0 (12):\penalty0 125017, 2017.

\bibitem[Hauptmann et~al.(2018{\natexlab{b}})Hauptmann, Lucka, Betcke, Huynh,
  Adler, Cox, Beard, Ourselin, and Arridge]{Hauptmann2018a}
Andreas Hauptmann, Felix Lucka, Marta Betcke, Nam Huynh, Jonas Adler, Ben Cox,
  Paul Beard, Sebastien Ourselin, and Simon Arridge.
\newblock Model-based learning for accelerated, limited-view 3-d photoacoustic
  tomography.
\newblock \emph{IEEE transactions on medical imaging}, 37\penalty0
  (6):\penalty0 1382--1393, 2018{\natexlab{b}}.

\bibitem[Hyv{\"o}nen(2009)]{Hyvoenen2009}
Nuutti Hyv{\"o}nen.
\newblock Approximating idealized boundary data of electric impedance
  tomography by electrode measurements.
\newblock \emph{Mathematical Models and Methods in Applied Sciences},
  19\penalty0 (07):\penalty0 1185--1202, 2009.

\bibitem[Isaacson et~al.(2004)Isaacson, Mueller, Newell, and
  Siltanen]{Isaacson2004}
D.~Isaacson, J.~L. Mueller, J.~C. Newell, and S.~Siltanen.
\newblock Reconstructions of chest phantoms by the {D}-bar method for
  electrical impedance tomography.
\newblock \emph{IEEE Transactions on Medical Imaging}, 23:\penalty0 821--828,
  2004.

\bibitem[Jin et~al.(2017)Jin, McCann, Froustey, and Unser]{Jin2017}
Kyong~Hwan Jin, Michael~T McCann, Emmanuel Froustey, and Michael Unser.
\newblock Deep convolutional neural network for inverse problems in imaging.
\newblock \emph{IEEE Transactions on Image Processing}, 26\penalty0
  (9):\penalty0 4509--4522, 2017.

\bibitem[Kang et~al.(2017)Kang, Min, and Ye]{Kang2017}
Eunhee Kang, Junhong Min, and Jong~Chul Ye.
\newblock A deep convolutional neural network using directional wavelets for
  low-dose x-ray ct reconstruction.
\newblock \emph{Medical physics}, 44\penalty0 (10), 2017.

\bibitem[Knudsen et~al.(2009)Knudsen, Lassas, Mueller, and
  Siltanen]{Knudsen2009}
K.~Knudsen, M.~Lassas, J.L. Mueller, and S.~Siltanen.
\newblock Regularized {D}-bar method for the inverse conductivity problem.
\newblock \emph{Inverse Problems and Imaging}, 3\penalty0 (4):\penalty0
  599--624, 2009.

\bibitem[Kourunen et~al.(2008)Kourunen, Savolainen, Lehikoinen, Vauhkonen, and
  Heikkinen]{Kourunen2008}
J~Kourunen, T~Savolainen, A~Lehikoinen, M~Vauhkonen, and LM~Heikkinen.
\newblock Suitability of a pxi platform for an electrical impedance tomography
  system.
\newblock \emph{Measurement Science and Technology}, 20\penalty0 (1):\penalty0
  015503, 2008.

\bibitem[Liu et~al.(2005)Liu, Saulnier, Newell, Isaacson, and Kao]{Liu2005}
N.~Liu, G.J. Saulnier, J.C. Newell, D.~Isaacson, and T-J. Kao.
\newblock Act4: a high-precision, multi-frequency electrical impedance
  tomograph.
\newblock Presented at 6th Conference on Biomedical Applications of Electrical
  Impedance Tomography, June 2005.
\newblock London, U.K.

\bibitem[Lytle et~al.(2018)Lytle, Perry, and Siltanen]{Lytle2018}
George Lytle, Peter Perry, and Samuli Siltanen.
\newblock Nachman's reconstruction method for the calder\'on problem with
  discontinuous conductivities.
\newblock \emph{(preprint) arXiv:1809.09272}, 2018.

\bibitem[Martin and Choi(2017)]{Martin2017}
S{\'e}bastien Martin and Charles~TM Choi.
\newblock A post-processing method for three-dimensional electrical impedance
  tomography.
\newblock \emph{Scientific reports}, 7\penalty0 (1):\penalty0 7212, 2017.

\bibitem[Mueller and Siltanen(2012)]{Mueller2012}
J.L. Mueller and S.~Siltanen.
\newblock \emph{Linear and Nonlinear Inverse Problems with Practical
  Applications}.
\newblock SIAM, 2012.

\bibitem[Mueller et~al.(2002)Mueller, Siltanen, and Isaacson]{Mueller2002}
J.L. Mueller, S.~Siltanen, and D.~Isaacson.
\newblock A direct reconstruction algorithm for electrical impedance
  tomography.
\newblock \emph{IEEE Transactions on Medical Imaging}, 21\penalty0
  (6):\penalty0 555--559, 2002.

\bibitem[Murphy and Mueller(2009)]{Murphy2009}
E.~K. Murphy and J.~L. Mueller.
\newblock Effect of domain-shape modeling and measurement errors on the 2-d
  {D}-bar method for electrical impedance tomography.
\newblock \emph{IEEE Transactions on Medical Imaging}, 28\penalty0
  (10):\penalty0 1576--1584, 2009.

\bibitem[Nachman(1996)]{Nachman1996}
A.~I. Nachman.
\newblock Global uniqueness for a two-dimensional inverse boundary value
  problem.
\newblock \emph{Annals of Mathematics}, 143:\penalty0 71--96, 1996.

\bibitem[Ronneberger et~al.(2015)Ronneberger, Fischer, and
  Brox]{Ronneberger2015}
Olaf Ronneberger, Philipp Fischer, and Thomas Brox.
\newblock U-net: Convolutional networks for biomedical image segmentation.
\newblock In \emph{International Conference on Medical image computing and
  computer-assisted intervention}, pages 234--241. Springer, 2015.

\bibitem[Rudin et~al.(1992)Rudin, Osher, and Fatemi]{Rudin1992}
L.I. Rudin, S.~Osher, and E.~Fatemi.
\newblock Nonlinear total variation based noise removal algorithms.
\newblock \emph{Physica D: Nonlinear Phenomena}, 60\penalty0 (1-4):\penalty0
  259--268, 1992.

\bibitem[Schlemper et~al.(2018)Schlemper, Caballero, Hajnal, Price, and
  Rueckert]{Schlemper2018}
Jo~Schlemper, Jose Caballero, Joseph~V Hajnal, Anthony~N Price, and Daniel
  Rueckert.
\newblock A deep cascade of convolutional neural networks for dynamic mr image
  reconstruction.
\newblock \emph{IEEE transactions on Medical Imaging}, 37\penalty0
  (2):\penalty0 491--503, 2018.

\bibitem[Seo et~al.(2018)Seo, Kim, Jargal, Lee, and Harrach]{Seo2018}
Jin~Keun Seo, Kang~Cheol Kim, Ariungerel Jargal, Kyounghun Lee, and Bastian
  Harrach.
\newblock A learning-based method for solving ill-posed nonlinear inverse
  problems: a simulation study of lung eit.
\newblock \emph{arXiv preprint arXiv:1810.10112}, 2018.

\bibitem[Siltanen and Tamminen(2016)]{Siltanen2014}
S.~Siltanen and J.~P. Tamminen.
\newblock Reconstructing conductivities with boundary corrected d-bar method.
\newblock \emph{Journal of Inverse and Ill-posed Problems}, 22\penalty0
  (6):\penalty0 847--870, 2016.

\bibitem[Somersalo et~al.(1992)Somersalo, Cheney, and Isaacson]{Somersalo1992}
Erkki Somersalo, Margaret Cheney, and David Isaacson.
\newblock Existence and uniqueness for electrode models for electric current
  computed tomography.
\newblock \emph{SIAM Journal on Applied Mathematics}, 52\penalty0 (4):\penalty0
  1023--1040, 1992.

\end{thebibliography}
\end{document}